\documentclass{amsart}
\usepackage{amscd}
\usepackage{amsfonts}
\usepackage{amsmath}
\usepackage{amssymb}
\usepackage{latexsym}
\usepackage{amsthm,color}
\usepackage[all]{xy}
\usepackage{epsfig}
\usepackage{graphicx}
\usepackage{color}
\usepackage{tikz}
\usepackage{bm}
\usepackage{hyperref}
\usetikzlibrary{cd}
\usetikzlibrary{calc}
\newtheorem{proposition}{Proposition}
\newtheorem{lemma}{Lemma}[section]

\newtheorem{theorem}{Theorem}
\newtheorem{corollary}{Corollary}

\theoremstyle{definition}
\newtheorem{definition}{Definition}
\theoremstyle{definition}
\newtheorem{example}{Example}
\theoremstyle{definition}

\newtheorem{note}{Note}

\newtheorem{problem}{Problem}

\begin{document}
\title[
Envelopes created by sphere families in Euclidean 3-space]
{Envelopes created by sphere families in Euclidean 3-space}
\author[T.~Nishimura]{Takashi Nishimura
}
\address{
Research Institute of Environment and Information Sciences,
Yokohama National University,
Yokohama 240-8501, Japan
}
\email{nishimura-takashi-yx@ynu.ac.jp}
\author[M.~Takahashi]{Masatomo Takahashi
}
\address{
Muroran Institute of Technology,
Muroran 050-8585, Japan
}
\email{masatomo@muroran-it.ac.jp}
\author[Y.~Wang]{Yongqiao Wang
}
\address{
School of Science,
Dalian Maritime University,
Dalian 116026,
P.R. China
}
\email{wangyq@dlmu.edu.cn}
\begin{abstract}
In this paper,  on envelopes created by sphere families in Euclidean 3-space,
all four basic problems (existence problem, representation problem,
problem on the number of envelopes, problem on relationships of definitions)
are solved.
\end{abstract}


\date{}

\maketitle

\section{Introduction\label{section1}}
Throughout this paper,
$U$ is
{\color{black} an open subset of $\mathbb{R}^2$}
and all functions,
mappings
{\color{black}and vector fields
are of class $C^\infty$,
}
unless otherwise stated.
\medskip

Since the dawn of differential analysis, the envelopes
{\color{black}created} by families of submanifolds have constituted a significant area of inquiry, with numerous applications in differential equations, differential geometry, engineering, and physics (see \cite{brucegiblin,craig,EFK,history, izumiya}). The study of envelopes of families of straight lines in the plane, which represents the classical case, dates back to Leibniz and Johann Bernoulli.
By using a geometric mechanism, the first author derived precise expressions for envelopes of straight line families and hyperplane families in a concise geometric way, without the need to solve the corresponding system of differential equations with a constraint condition (see \cite{nishimura,nishimura1}).
Similarly to straight lines, circles are well-known plane curves, and envelopes of families of circles have also been investigated (see \cite{WN}).

On the other hand, the study of sphere families in Euclidean 3-space should not be overlooked, particularly because their envelopes play a crucial role in applications such as {\color{black}the} seismic survey. The following is a brief explanation of this technique, based on 7.14(9) of \cite{brucegiblin}.
In Euclidean 3-space, let the ``ground level'' $G$ be parameterized by $\bm{x}:U\rightarrow\mathbb{R}^3$. Suppose a granite stratum lies beneath an upper layer of sandstone, with their interface $M$ parameterized by $\widetilde{\bm{f}}:U\rightarrow\mathbb{R}^3$.
{\color{black}The s}eismic survey is a method for approximating $M$ as accurately as possible. The procedure involves detonating a source at a fixed point $\bm{P}$ on $G$. Under the assumption that sound waves propagate linearly, they reflect off $M$ and return to sensors positioned along $G$. These sensors, located at points $\bm{x}(u,v)$ {\color{black}$((u, v)\in U)$}, precisely record the arrival times of the reflected waves (see Fig. 1).
\begin{figure}[h]
\begin{center}
\includegraphics[width=10cm]
{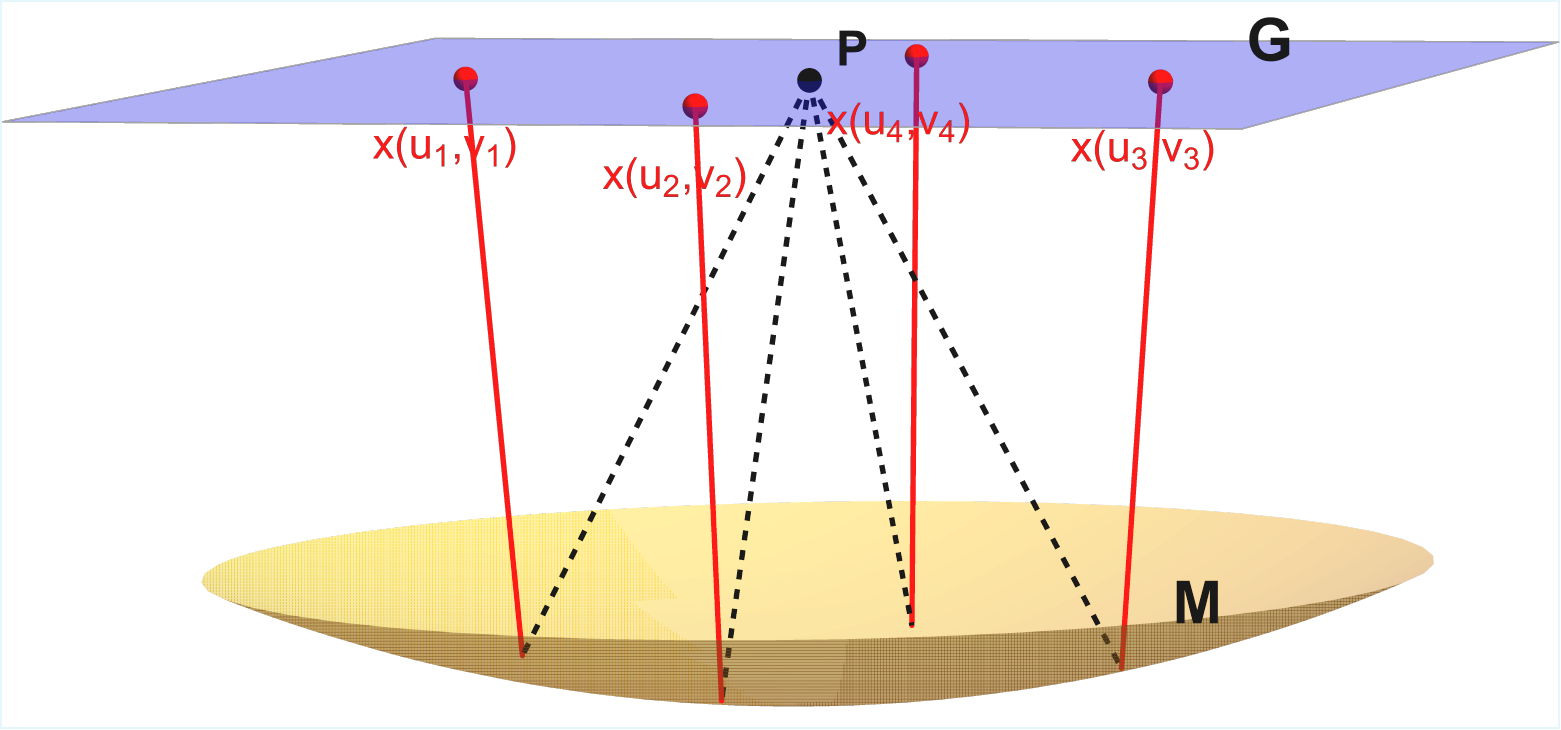}
\caption{Reflection of sound waves.
}
\label{seismic survey 1}
\end{center}
\end{figure}

It is known that there exists a surface $W$, parameterized by $\bm{f}:U\rightarrow\mathbb{R}^3$
with well-defined normals, such that
any broken line of a reflected ray path{\color{black}ing} from
{\color{black}the} point $\bm{P}$ to the surface $M$ can be replaced by a straight line segment of equal length that is normal to $W$. The surface $W$ is termed the {\it orthotomic} of $M$ relative to $\bm{P}$, and conversely, $M$ is called the {\it anti-orthotomic} of $W$ relative to $\bm{P}$. This implies that $W$ can be geometrically reconstructed as the envelope of the family of spheres centered at points $\bm{x}(u,v)$ on $G$ with radii equal to the distances $\|\bm{f}(u,v)-\bm{x}(u,v)\|$ (see Fig. 2).
\begin{figure}[h]
\begin{center}
\includegraphics[width=11cm]
{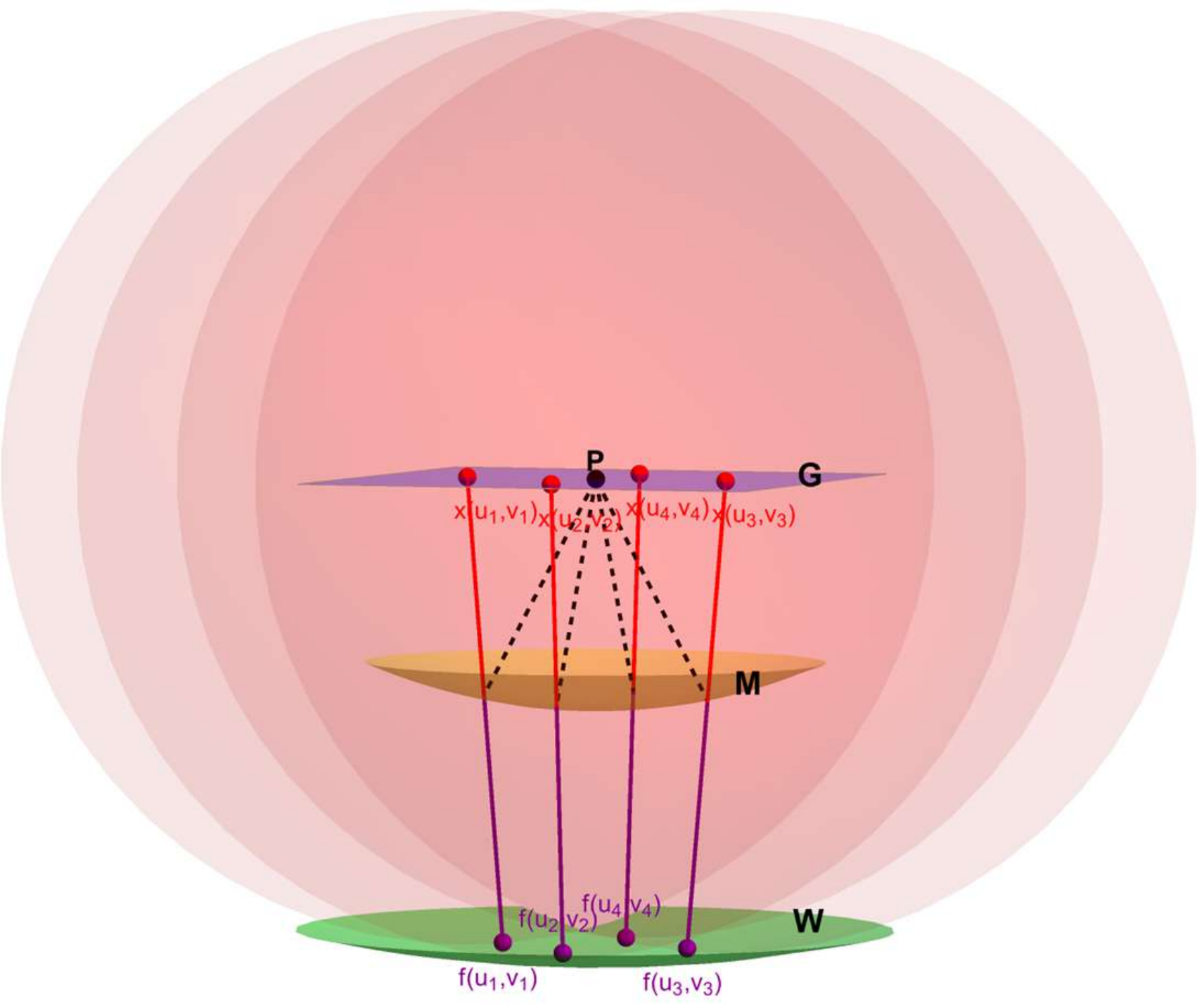}
\caption{An envelope created by the sphere family.
}
\label{seismic survey 2}
\end{center}
\end{figure}
Once the parametrization $\bm{f}$ of $W$ is obtained, the parametrization of $M$ can then be readily derived using the anti-orthotomic technique established in \cite{STN}. Therefore, in order to investigate the parametrization of $W$ as precisely as
 possible, it is very important to construct a general theory
 on envelopes created by sphere families in Euclidean 3-space, which is the main purpose of this paper.
{\color{black}More precisely, in the main theorem (Theorem \ref{theorem1}),
given a family of spheres in $\mathbb{R}^3$, it is provided
(1) a necessary and sufficient condition for the existence of its envelope,
(2) a representation formula of the envelope in terms of the given
family of spheres when
an envelope exists and
(3) the determination of \lq\lq number\rq\rq\; of envelopes
when an envelope exists.}
{\color{black}Moreover, in the second main theorem
(Theorem \ref{theorem3}), it is provided a description
on the explicit difference of two apparently different definitions of envelope for a given family of spheres creating an envelope.
}
\medskip

The rest of the paper is organized as follows: First, we introduce some necessary definitions in Section 2, including the definition of framed surfaces, which was originally provided by the second author and his collaborators. The main result Theorem \ref{theorem1} is also presented in this section. Then we give the proof of Theorem \ref{theorem1} in Section 3. In Section 4, some examples are given to show how
Theorem \ref{theorem1} is effectively applicable.
In Section 5, we investigate the
{\color{black}
relations}
of
{\color{black}
two distinct definitions}
of an envelope created by a sphere family.
{\color{black}
The second main result (Theorem \ref{theorem3}) is presented
in this section.
}
Finally
we apply Theorem \ref{theorem1}
to study the associated surfaces of framed surfaces
{\color{black}in Section 6}.

\section{Basic concepts and the main result}
\label{sectiontwo}   %
For a point $\bm{Q}$ of $\mathbb{R}^3$
and a positive number $\lambda$,
the sphere $S_{(Q, \lambda)}$ with centre $\bm{Q}$ and radius $\lambda$
is naturally defined as follows.
\[
S_{(Q, \lambda)}=
\left\{\left.
(x_1, x_2, x_3)\in \mathbb{R}^3\, \right|\,
\left((x_1, x_2, x_3) - \bm{Q}\right)\cdot \left((x_1, x_2, x_3) - \bm{Q}\right) = \lambda^2
\right\},
\]
where
the dot in the center
stands for the standard scalar product of $3$-dimensional vectors.
For a surface $\bm{x}: U\to \mathbb{R}^3$ and a positive function
$\lambda: U\to \mathbb{R}_+$,
the sphere family $\mathcal{S}_{(x, \lambda)}$
is naturally defined as follows.
Here, $\mathbb{R}_+$ stands for the set consisting of positive real numbers.
\[
\mathcal{S}_{(x, \lambda)}=
\left\{S_{(x(u,v), \lambda(u,v))}
\right\}_{(u,v)\in U}.
\]
Let $S^2=\{\bm{a}\in\mathbb{R}^3|\bm{a}\cdot \bm{a}=1\}$ be the unit sphere in $\mathbb{R}^3$. We denote the $3$-dimensional smooth manifold
$\{(\bm{a}{\color{black},\,}\bm{b})\in S^2\times S^2|\bm{a}\cdot \bm{b}=0\}$ by $\Delta$. The definitions of framed surfaces and framed base surfaces given in \cite{fukunagatakahashi1} are as follows.
\begin{definition}\label{framed surface}(\cite{fukunagatakahashi1})
{\rm
We call $(\bm{x},\bm{n},\bm{s}):U\rightarrow\mathbb{R}^3\times\Delta$ a \textit{framed surface} if
$\bm{x}_u(u,v)\cdot \bm{n}(u,v)=0$ and $\bm{x}_v(u,v)\cdot \bm{n}(u,v)=0$ hold for all $(u,v)\in U,$ where $\bm{x}_u(u,v)=(\partial \bm{x}/\partial u)(u,v)$ and $\bm{x}_v(u,v)=(\partial \bm{x}/\partial v)(u,v)$. We say that $\bm{x}:U\rightarrow\mathbb{R}^3$ is a \textit{framed base surface} if there exists a mapping $(\bm{n},\bm{s}):U\rightarrow\Delta$ such that $(\bm{x},\bm{n},\bm{s})$ is a framed surface.
}
\end{definition}
{\color{black}
\noindent
On the other hand, the notion of
frontal is defined in a manner analogous to that of a framed
{\color{black}base}
surface, as follows.
\begin{definition}\label{frontal}
We say that $\bm{x}:U\rightarrow\mathbb{R}^3$ is a
\textit{frontal} if there exists a mapping
$\bm{n}:U\rightarrow S^2$
such that
$\bm{x}_u(u,v)\cdot \bm{n}(u,v)=0$ and $\bm{x}_v(u,v)\cdot \bm{n}(u,v)=0$
hold for all $(u,v)\in U$.
\end{definition}
}
\par
\noindent
For frontals,
{\color{black}
\cite{ishikawa, ishikawa2}
}
are recommended as excellent references.
{\color{black}
\begin{proposition}\label{proposition0}
\label{proposition_equivalence}
{\color{black}
Any framed base surface
$\bm{x}:U\rightarrow\mathbb{R}^3$ is a frontal.
Conversely, suppose that $\bm{x}: U\to \mathbb{R}^3$ is a frontal.   Namely,
suppose that there exists a mapping
$\bm{n}:U\rightarrow S^2$
such that
$\bm{x}_u(u,v)\cdot \bm{n}(u,v)=0$ and $\bm{x}_v(u,v)\cdot \bm{n}(u,v)=0$
hold for all $(u,v)\in U$.
Then, one of the following two holds.
\begin{enumerate}
\item[(1)] The frontal $\bm{x}$ is a framed base surface.
\item[(2)] For any point $\bm{X}_0\in S^2$,
the restriction of $\bm{x}$ to $U\setminus\bm{n}^{-1}(\bm{X}_0)$
is a framed base surface.
\end{enumerate}
}
\end{proposition}
\begin{proof}
{\color{black}
By definition, any framed base surface
$\bm{x}:U\rightarrow\mathbb{R}^3$ is clearly a frontal.
\par
Next, suppose that  $\bm{x}:U\rightarrow\mathbb{R}^3$ is a frontal and
the assertion (1) does not hold.
Let $\bm{X}_0$ be a point of $S^2$.
Set
\[
U_{\bm{X}_0}= U\setminus\bm{n}^{-1}(\bm{X}_0), \:
\bm{x}_{\bm{X}_0}=\bm{x}|_{U_{\bm{X}_0}}\, \mbox{ and }\,
 \bm{n}_{\bm{X}_0}=\bm{n}|_{U_{\bm{X}_0}}.
\]
Then, it is clear that
\[
\bm{n}_{\bm{X}_0}\left(U_{\bm{X}_0}\right)\subset S^2\setminus\bm{X}_0.
\]
Let $\Phi_{\bm{X}_0}: S^2\setminus\bm{X}_0 \to \Pi_{\bm{X}_0}$
be the stereographic projection, where
$\Pi_{\bm{X}_0}$
{\color{black}is the
plane}
in
the $3$-dimensional
{\color{black}affine space}
$\mathbb{R}^3$ defined by
\[
\Pi_{\bm{X}_0} = \left\{\left.\bm{X}\in \mathbb{R}^3\, \right|\,
\bm{X}\cdot
{\color{black}
\left(-\bm{X}_0\right)
}
=1\right\}.
\]
{\color{black}
The plane $\Pi_{\bm{X}_0}$ may be canonically identified with the $2$-dimensional
vector space $\mathbb{R}^2$.
}
Let $\left\{\bm{e}_1,\, \bm{e}_2\right\}$ be a basis of $\Pi_{\bm{X}_0}$.
Since $\Pi_{\bm{X}_0}$ is an affine
{\color{black}plane, it follows that}
$\left\{\bm{e}_1,\, \bm{e}_2\right\}$
{\color{black}
may be regarded as}
a basis of $T_{\bm{X}}\Pi_{\bm{X}_0}$
as well for any $\bm{X}\in \Pi_{\bm{X}_0}$.
Since the stereographic projection $\Phi_{\bm{X}_0}$ is a diffeomorphism,
$d\left(\Phi_{\bm{X}_0}\right)
{\color{black}_{\widetilde{\bm{X}}}^{-1}}
\left(\bm{e}_i\right)$ is a non-zero
vector in $T_{\widetilde{\bm{X}}}S^2$ for any $\widetilde{\bm{X}}\in S^2
\setminus\bm{X}_0
$ and any $i\in \{1,2\}$.
Moreover, from the construction it follows
\[
\widetilde{\bm{X}}\cdot d\left(\Phi_{\bm{X}_0}\right)_{\widetilde{\bm{X}}}^{-1}\left(\bm{e}_i\right)=0
\]
for any $\widetilde{\bm{X}}\in S^2
\setminus\bm{X}_0
$ and any $i\in \{1,2\}$.
For any $
{\color{black}
(u,v)
}
\in U_{\bm{X}_0}$,
set
\[
\bm{s}_{\bm{X}_0}{\color{black}(u,v)} =
\frac{d\left(\Phi_{\bm{X}_0}\right)_{\bm{n}_{\bm{X}_0}
{\color{black}(u,v)}}^{-1}\left(\bm{e}_1\right)}
{||d\left(\Phi_{\bm{X}_0}\right)_{\bm{n}_{\bm{X}_0}
{\color{black}(u,v)}}^{-1}
\left(\bm{e}_1\right)||},
\]
where $||\bm{a}||$ denotes the length of the vector $\bm{a}\in T_{\bm{n}_{\bm{X}_0}(u,v)}S^2$.
Then, it is clear that $\bm{s}_{\bm{X}_0}
(u,v)
\in S^2$.
Therefore, $\left(\bm{x}_{\bm{X}_0}, \bm{n}_{\bm{X}_0}, \bm{s}_{\bm{X}_0}\right)$ is a framed surface and $\bm{x}_{\bm{X}_0}$
is a framed base surface.
}
\end{proof}
}
\par
{\color{black}
\begin{note}\label{note0}
\item[(1)] Notice that in Proposition \ref{proposition0},
if the assertion (1) holds, then the assertion (2), too, clearly holds.
Thus, by Proposition \ref{proposition0}, the assertion (2) always holds.   \\
\item[(2)]
For any $\bm{X}_0\in S^2$, consider the exponential mapping
\[
\exp_{\bm{X}_0}: \Pi_{\bm{X}_0}\to S^2.
\]
Set
\[
U  =  \Pi_{\bm{X}_0}, \;
\bm{x}   =  \exp_{\bm{X}_0}.
\]
Then, the mapping $\bm{x}: U\to S^2\subset \mathbb{R}^3$
is a frontal satisfying $\bm{x}\left(U\right)=S^2$.
It is well-known that there are no
non-zero tangent vector fields over $S^2$.   Thus,
for the frontal $\bm{x}: U\to S^2\subset \mathbb{R}^3$,
the assertion (1) of Proposition \ref{proposition0} does not hold.
On the other hand, by Proposition \ref{proposition0},
the assertion (2) holds for the same frontal $\bm{x}$.
Therefore, as for Proposition \ref{proposition0},
the assertion (1) is certainly stronger than the assertion (2).
\end{note}
}

Let $(\bm{x},\bm{n},\bm{s}):U\rightarrow\mathbb{R}^3\times\Delta$ be a framed surface, we denote the vector product of $\bm{n}(u,v)$ and $\bm{s}(u,v)$ by $\bm{t}(u,v)=\bm{n}(u,v)\times \bm{s}(u,v)$. Then $\{\bm{n}(u,v),\bm{s}(u,v),\bm{t}(u,v)\}$ is a moving frame along $\bm{x}(u,v)$, and we have the following systems of differential equations:
\begin{align}
\left(
  \begin{array}{c}
    \bm{x}_u \\
    \bm{x}_v \\
  \end{array}
\right)&=\left(
          \begin{array}{cc}
            a_1 & b_1 \\
            a_2 & b_2 \\
          \end{array}
        \right)
        \left(
          \begin{array}{c}
            \bm{s} \\
            \bm{t} \\
          \end{array}
        \right),
\\
\left(
  \begin{array}{c}
    \bm{n}_u \\
    \bm{s}_u \\
    \bm{t}_u \\
  \end{array}
\right)&=\left(
          \begin{array}{ccc}
            0 & e_1 & f_1 \\
            -e_1 & 0 & g_1 \\
            -f_1 & -g_1 & 0  \\
          \end{array}
        \right)
        \left(
          \begin{array}{c}
            \bm{n}\\
            \bm{s} \\
            \bm{t}\\
          \end{array}
        \right),\\
\left(
  \begin{array}{c}
    \bm{n}_v \\
    \bm{s}_v \\
    \bm{t}_v \\
  \end{array}
\right)&=\left(
          \begin{array}{ccc}
            0 & e_2 & f_2 \\
            -e_2 & 0 & g_2 \\
            -f_2 & -g_2 & 0  \\
          \end{array}
        \right)
        \left(
          \begin{array}{c}
            \bm{n}\\
            \bm{s} \\
            \bm{t} \\
          \end{array}
        \right),
\end{align}
where $a_i$, $b_i$, $e_i$, $f_i$, $g_i$ ($i=1,2$) are
{\color{black}$C^\infty$}
functions and we call these functions {\it basic invariants} of the framed surface. We denote the above three matrices by $\mathcal{A}$, $\mathcal{F}_1$ and $\mathcal{F}_2$ and also call them {\it basic invariants} of the framed surface.
{\color{black}Notice that}
$(u_0,v_0)
{\color{black}\in U}$
is a singular point of $\bm{x}$ if and only if $\mathrm{det}\mathcal{A}=0.$ According to the integrability conditions $\bm{x}_{uv}=\bm{x}_{vu}$ and $\mathcal{F}_{2,u}-\mathcal{F}_{1,v}=\mathcal{F}_1\mathcal{F}_2-\mathcal{F}_2\mathcal{F}_1$, we have the following equalities for basic invariants:
\begin{align}
\left\{\begin{array}{c}
        a_{1,v}-b_1g_2=a_{2,u}-b_2g_1, \\
        b_{1,v}-a_2g_1=b_{2,u}-a_1g_2, \\
        a_1e_2+b_1f_2=a_2e_1+b_2f_1,
      \end{array}\right.
\end{align}
\begin{align}
\left\{\begin{array}{c}
        e_{1,v}-f_1g_2=e_{2,u}-f_2g_1, \\
        f_{1,v}-e_2g_1=f_{2,u}-e_1g_2, \\
        g_{1,v}-e_1f_2=g_{2,u}-e_2f_1.
      \end{array}\right.
\end{align}
\medskip
\par
In this paper,
{\color{black}
to be consistent with \cite{WN},
}
the surface $\bm{x}:U\to \mathbb{R}^3$
for a sphere family $\mathcal{S}_{(x, \lambda)}$ is assumed to be a
{\color{black}frontal.
By the assertion (2) of Proposition \ref{proposition_equivalence},
we may say also that  {\color{black}for any $\bm{X}_0\in S^2$} the surface
$\bm{x}_{\bm{X}_0}:U_{\bm{X}_0}\to \mathbb{R}^3$
for a sphere family $\mathcal{S}_{(x, \lambda)}$
is assumed to be a framed base surface.
{\color{black}Notice that there exists a point
$\bm{X}_0\in S^2$ such that $\bm{n}^{-1}(\bm{X}_0)=\emptyset$
when we apply the results provided in this paper to the seismic survey.
Thus we may say that
the assertion (1) of Proposition
\ref{proposition_equivalence} holds in the case of application
to the seismic survey.
In addition, when we apply the results of this paper
to real-world problems, the parameter spaces $U$ are
approximately constructed from discrete data.
Hence, we may say that in the case of practical real-world applications,
the surface $\bm{x}:U\to \mathbb{R}^3$
 for a sphere family $\mathcal{S}_{(x, \lambda)}$ is assumed to be
 a framed base surface.}}
The following is adopted as the definition of an envelope created by a sphere family.
\begin{definition}\label{envelope}
{\rm
Let $\mathcal{S}_{(x, \lambda)}$ be a sphere family.
A mapping $\bm{f}: U\to \mathbb{R}^3$ is called an
\textit{envelope} created by $\mathcal{S}_{(x, \lambda)}$
if the following three hold for any $(u,v)\in U$.
\begin{enumerate}
\item[(1)] $\bm{f}(u,v)\in S_{(x(u,v), \lambda(u,v))}$.
\item[(2)] $\bm{f}_u(u,v)\cdot\left(\bm{f}(u,v)-\bm{x}(u,v)\right)=0$.
\item[(3)] $\bm{f}_v(u,v)\cdot\left(\bm{f}(u,v)-\bm{x}(u,v)\right)=0$.
\end{enumerate}
}
\end{definition}
\noindent
By definition, as similar as the envelope created by a hyperplane family (see \cite{nishimura}), the envelope created by a sphere family is a mapping giving a solution of two first order differential equations with one constraint condition. Moreover, from the definition again, it can be seen that $(\bm{f},\frac{\bm{f}-\bm{x}}{\|\bm{f}-\bm{x}\|}):U\rightarrow\mathbb{R}^3\times S^2$ is a Legendre surface, i.e., $\bm{f}:U\rightarrow\mathbb{R}^3$ is a frontal.

\begin{problem}\label{problem}
Let $(\bm{x},\bm{n},\bm{s}):U\rightarrow\mathbb{R}^3\times\Delta$ be a framed surface and let $\lambda:U\to \mathbb{R}_+$ be a positive function.
\begin{enumerate}
\item[(1)]
Find a necessary and sufficient condition for
the sphere family $\mathcal{S}_{(x, \lambda)}$ to create an envelope
in terms of $\bm{x}$, $\bm{n}$, $\bm{s}$ and $\lambda$.
\item[(2)] Suppose that the sphere family $\mathcal{S}_{(x, \lambda)}$
creates an envelope.    Then,
find a parametrization of the envelope
created by $\mathcal{S}_{(x, \lambda)}$
in terms of $\bm{x}$, $\bm{n}$, $\bm{s}$ and $\lambda$.
\item[(3)] Suppose that the sphere family $\mathcal{S}_{(x, \lambda)}$
creates an envelope.    Then, characterize the number of distinct envelopes
created by $\mathcal{S}_{(x, \lambda)}$
in terms of $\bm{x}$, $\bm{n}$, $\bm{s}$ and $\lambda$.
\end{enumerate}
\end{problem}

\begin{note}\label{note1}
\item[(1)] The item (1) of Problem \ref{problem} is a problem to seek the
integrability conditions.  There are various cases, for instance
the concentric sphere family
$$
\{\{(x_1, x_2, x_3)\in \mathbb{R}^3\, |\, x_1^2+x_2^2+x_3^2=u^2+v^2\}\}_{(u,v)
\in \mathbb{R}^2
{\color{black}
\setminus\{(0,0)\}}
}
$$
does not create an envelope while the parallel-translated sphere family
$$\{\{(x_1, x_2, x_3)\in \mathbb{R}^3\, |\, (x_1-u)^2+(x_2-v)^2+x_3^2=1\}\}_{(u,v)\in \mathbb{R}^2}$$
does create two envelopes.   Thus, (1) of Problem \ref{problem} is significant.
\item[(2)] The following Example \ref{example1} shows that
the well-known method to represent the envelope is useless in this case.
Thus, (2) of Problem \ref{problem} is important and the positive answer to it is
much desired.
\item[(3)]  The following Example \ref{example2} shows that
there are at least three cases: the case having a unique envelope, the case
having exactly two envelopes and the case having uncountably many envelopes.
Thus, (3) of Problem \ref{problem} is interesting and meaningful.
\end{note}

\begin{example}\label{example1}
Let $\bm{x}:\mathbb{R}^2\to \mathbb{R}^3$ be the mapping defined by
$\bm{x}(u,v)=\left(u^2, u^3, v\right)$. Set
$(\bm{n},\bm{s}):\mathbb{R}^2\rightarrow\Delta$, $\bm{n}(u,v)=(1/\sqrt{9u^2+4})(-3u,2,0)$, $\bm{s}=(0,0,1)$, then $(\bm{x},\bm{n},\bm{s}):\mathbb{R}^2\rightarrow\mathbb{R}^3\times\Delta$ is a framed surface.
Let $\lambda:\mathbb{R}^2\to \mathbb{R}_+$ be the constant function
defined by $\lambda(u,v)=1$.
Then, it seems that the sphere family
$\mathcal{S}_{(x, \lambda)}$ creates envelopes.
Thus, we can expect that the created envelopes can be obtained
by the well-known method.
Set $F(x_1, x_2, x_3, u, v)=\left(x_1-u^2\right)^2+\left(x_2-u^3\right)^2+\left(x_3-v\right)^2-1$.
Then, we have the following.
{\small
\begin{align*}
\mathcal{D}=&
\bigg\{(x_1, x_2, x_3)\in \mathbb{R}^3\, \left|\,
\exists u,v \mbox{ such that }F(x_1, x_2, x_3, u, v)=\frac{\partial F}{\partial u}(x_1, x_2, x_3, u, v)=\frac{\partial F}{\partial v}(x_1, x_2, x_3, u, v)=0
\right.\bigg\} \\
=&
\bigg\{(x_1, x_2, x_3)\in \mathbb{R}^3\, \left|
\exists u,v \mbox{ such that }
\left(x_1-u^2\right)^2+\left(x_2-u^3\right)^2+\left(x_3-v\right)^2-1=0, \, x_3-v=0, \right.\\
&
2u\left(x_1-u^2\right)+3u^2\left(x_2-u^3\right)=0
\bigg\} \\
=&
\left\{(x_1, x_2, x_3)\in \mathbb{R}^3\, \left|\,
x_1^2+x_2^2=1,\, x_3=v\right.\right\} \bigcup
\bigg\{(x_1, x_2, x_3)\in \mathbb{R}^3\, \left|\,
\left(x_1-u^2\right)^2+\left(x_2-u^3\right)^2-1=0,\right. \\
&
2(x_1-u^2)+3u(x_2-u^3)=0, \,
x_3=v
\bigg\} \\
= &
\left\{(x_1, x_2, x_3)\in \mathbb{R}^3\, \left|\,
x_1^2+x_2^2=1\right.\right\} \bigcup
\left\{\left.
\left(u^2\pm \frac{3u}{\sqrt{9u^2+4}}, u^3\mp \frac{2}{\sqrt{9u^2+4}}, v\right)
\in \mathbb{R}^3\, \right|\, (u,v)\in \mathbb{R}^2
\right\}.
\end{align*}
}
In Example \ref{example3} of Section \ref{section3},
it turns out that the set $\mathcal{D}$ calculated here is actually larger than
the set of envelopes created by $\mathcal{S}_{(x, \lambda)}$,
namely the cylinder $\left\{(x_1, x_2, x_3)\in \mathbb{R}^3\, \left|\,
x_1^2+x_2^2=1\right.\right\}$ is redundant.
Therefore,
the well-known method
to represent the envelopes does not
work well in general and it must be applied under
appropriate assumptions even for sphere families.
The sphere family $\mathcal{S}_{(x, \lambda)}$
and the candidates of its envelope are depicted in
Figure \ref{figure_example1}.
\begin{figure}[h]
\begin{center}
\includegraphics[width=7cm]
{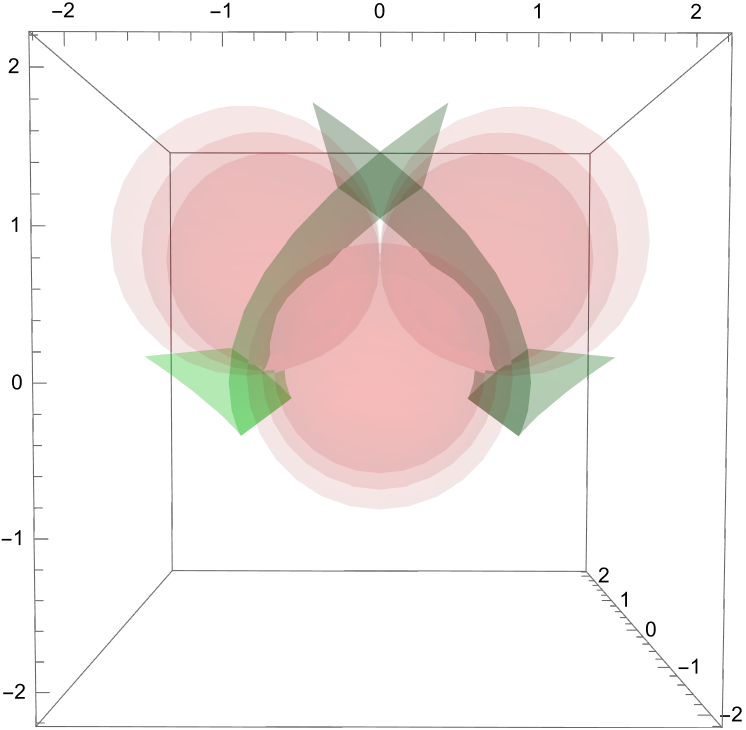}
\caption{The sphere family $\mathcal{S}_{(x, \lambda)}$
and the candidates of its envelope.
}
\label{figure_example1}
\end{center}
\end{figure}
\end{example}

\begin{example}\label{example2}
\begin{enumerate}
\item[(1)] Let $\bm{x}: \mathbb{R}^2
{\color{black}
\setminus\{(0,0)\}
}
\to \mathbb{R}^3$
be the mapping defined by
$\bm{x}(u,v)=\left(u, v, 0\right)$.  Then, it is clear that $\bm{x}$ is a framed base surface.
Let $\lambda: \mathbb{R}^2
{\color{black}
\setminus\{(0,0)\}
}
\to \mathbb{R}_+$ be the positive function
defined by $\lambda(u,v)=\sqrt{u^2+v^2}$.
Then, it is easily seen that the origin $(0,0,0)$ is an envelope created by the sphere family
$\mathcal{S}_{(x, \lambda)}$
and that there are no other envelopes created by
$\mathcal{S}_{(x, \lambda)}$.
Hence, the number of created envelopes is one
in this case.
\item[(2)] The parallel-translated sphere family
$\{\{(x_1, x_2, x_3)\in \mathbb{R}^3\, |\, (x_1-u)^2+(x_2-v)^2+x_3^2=1\}\}_{(u,v)\in \mathbb{R}^2}$
creates exactly two envelopes.
\item[(3)]   Let $\bm{x}: \mathbb{R}^2\to \mathbb{R}^3$ be the constant
mapping defined by $\bm{x}(u,v)=\left(0, 0, 0\right)$.
Then, it is clear that $\bm{x}$ is a framed base surface.
Let $\lambda: \mathbb{R}^2\to \mathbb{R}_+$ be the constant function
defined by $\lambda(u,v)=1$.
Then, for any two functions $\theta: \mathbb{R}^2\to \mathbb{R}$ and $\varphi: \mathbb{R}^2\to \mathbb{R}$,
the mapping $\bm{f}: \mathbb{R}^2\to \mathbb{R}^3$ defined by
$\bm{f}(u,v)=(\sin\varphi(u,v),\cos\varphi(u,v)\cos\theta(u,v), \sin\varphi(u,v)\cos\theta(u,v))$ is an envelope created by
the sphere family $\mathcal{S}_{(x, \lambda)}$.
Hence, there  are uncountably many created envelopes in this case.
\end{enumerate}
\end{example}

In order to solve Problem \ref{problem}, we define
the key definition of this paper as follows.
\begin{definition}\label{creative}
{\rm
Let $(\bm{x},\bm{n},\bm{s}): U\to \mathbb{R}^3\times\Delta$,  $\lambda: U\to \mathbb{R}_+$
be a framed surface and a positive function respectively.
Then, the sphere family $\mathcal{S}_{(x, \lambda)}$ is said to be
\textit{creative} if there exists a
mapping $\bm{\nu}:U\to S^2$ such that the following identities
hold for any $(u,v)\in U$.
\begin{enumerate}
\item[(1)] $\bm{x}_u(u,v)\cdot\bm{\nu}(u,v)+\lambda_u(u,v)=0$.
\item[(2)] $\bm{x}_v(u,v)\cdot\bm{\nu}(u,v)+\lambda_v(u,v)=0$.
\end{enumerate}
Set $\bm{\nu}(u,v)=\alpha(u,v)\bm{s}(u,v)+\beta(u,v)\bm{t}(u,v)+\gamma(u,v)\bm{n}(u,v)$. Then, the creative condition is equivalent to say that there exists a mapping $(\alpha,\beta,\gamma):U\to\mathbb{R}^3$ such that the following three identities hold for any $(u,v)\in U$.
\begin{enumerate}
\item[(1)] $a_1(u,v)\alpha(u,v)+b_1(u,v)\beta(u,v)+\lambda_u(u,v)=0$.
\item[(2)] $a_2(u,v)\alpha(u,v)+b_2(u,v)\beta(u,v)+\lambda_v(u,v)=0$.
\item[(3)] $\alpha^2(u,v)+\beta^2(u,v)+\gamma^2(u,v)=1$.
\end{enumerate}
}
\end{definition}

\begin{proposition}
Suppose the sphere family $\mathcal{S}_{(x, \lambda)}$ is creative. Let $(\alpha,\beta,\gamma):U\to\mathbb{R}^3$ be a mapping satisfying the creative condition, then we have
{\color{black}
\begin{eqnarray*}
\alpha(u,v)\mathrm{det}\left(\begin{array}{cc}
  b_1(u,v) & g_1(u,v) \\
  b_2(u,v) & g_2(u,v)
\end{array}\right)
& + &
\beta(u,v)\mathrm{det}\left(\begin{array}{cc}
  g_1(u,v) & a_1(u,v) \\
  g_2(u,v) & a_2(u,v)
\end{array}\right) \\
{ } & + &
\mathrm{det}\left(\begin{array}{cc}
  a_1(u,v) & \alpha_u(u,v) \\
  a_2(u,v) & \alpha_v(u,v)
\end{array}\right)+\mathrm{det}\left(\begin{array}{cc}
  b_1(u,v) & \beta_u(u,v) \\
  b_2(u,v) & \beta_v(u,v)
\end{array}\right)=0.
\end{eqnarray*}
}
\end{proposition}
\begin{proof}
By the creative condition, we have
\begin{align*}
a_1\alpha+b_1\beta+\lambda_u=0,~~
a_2\alpha+b_2\beta+\lambda_v=0.
\end{align*}
It follows
\begin{align*}
(a_{1,v}-a_{2,u})\alpha+(b_{1,v}-b_{2,u})\beta=a_2\alpha_u-a_1\alpha_v+b_2\beta_u-b_1\beta_v
\end{align*}
from $\lambda_{uv}=\lambda_{vu}$. Moreover, according to the integrability conditions $a_{1,v}-b_1g_2=a_{2,u}-b_2g_1$ and $b_{1,v}-a_2g_1=b_{2,u}-a_1g_2$, then the conclusion holds.
\end{proof}

We consider a parameter change of the domain $U$. Let $(\bm{x},\bm{n},\bm{s}): U\to \mathbb{R}^3\times\Delta$ be a framed surface with the basic invariants $\mathcal{A}$, $\mathcal{F}_1$ and $\mathcal{F}_2$.
{\color{black}
Let $V$ be an open subset of $\mathbb{R}^2$
which is $C^\infty$ diffeomorphic to
$U$.}
Suppose that $\phi:V\rightarrow U$, $(p,q)\mapsto\phi(p,q)=(u(p,q),v(p,q))$ is a parameter change, $(\widetilde{\bm{x}},\widetilde{\bm{n}},\widetilde{\bm{s}})=(\bm{x},\bm{n},\bm{s})\circ\phi: V\to \mathbb{R}^3\times\Delta$ has been proven to be a framed surface in \cite{fukunagatakahashi1}. Moreover, we have the following lemma.
\begin{lemma}\label{parameter le} (\cite{fukunagatakahashi1})
Under the above notations, the basic invariants
of $(\widetilde{\bm{x}},\widetilde{\bm{n}},\widetilde{\bm{s}})$ is given by
\begin{align*}
\left(\begin{array}{cc}
  \widetilde{a_1} & \widetilde{b_1} \\
  \widetilde{a_2} & \widetilde{b_2}
\end{array}\right)(p,q)=&\left(\begin{array}{cc}
  u_p & v_p \\
  u_q & v_q
\end{array}\right)(p,q)\left(\begin{array}{cc}
  a_1 & b_1 \\
  a_2 & b_2
\end{array}\right)(\phi(p,q)),\\
\left(\begin{array}{ccc}
        \widetilde{e_1} & \widetilde{f_1} & \widetilde{g_1} \\
        \widetilde{e_2} & \widetilde{e_2} & \widetilde{g_2}
      \end{array}
\right)(p,q)=&\left(\begin{array}{cc}
  u_p & v_p \\
  u_q & v_q
\end{array}\right)(p,q)\left(\begin{array}{ccc}
        e_1 & f_1 & g_1 \\
        e_2 & e_2 & g_2
      \end{array}
\right)(\phi(p,q)).
\end{align*}
\end{lemma}

\begin{proposition}
Under the above notations,
the sphere family $\mathcal{S}_{(x, \lambda)}$ is creative if and only if the sphere family $\mathcal{S}_{(\widetilde{x}, \widetilde{\lambda})}$ is creative, where $\widetilde{\lambda}(p,q)=\lambda\circ\phi(p,q)$.
\end{proposition}
\begin{proof}
By Definition \ref{creative}, it follows that $\mathcal{S}_{(x, \lambda)}$ is creative if and only if there exists a mapping $(\alpha,\beta,\gamma):U\to S^2$ such that
\begin{align*}
\left(\begin{array}{cc}
        a_1 & b_1 \\
        a_2 & b_2
      \end{array}
\right)(\phi(p,q))\left(\begin{array}{c}
                          \alpha \\
                          \beta
                        \end{array}
\right)(\phi(p,q))+\left(\begin{array}{c}
                          \lambda_u \\
                          \lambda_v
                        \end{array}
\right)(\phi(p,q))=\left(\begin{array}{c}
                          0 \\
                          0
                        \end{array}
\right).
\end{align*}
Because $\phi:V\rightarrow U$ is a parameter change, we have $\mathrm{det}(\phi_p(p,q),\phi_q(p,q))\neq0$ for any $(p,q)\in V$. It follows that the above equation holds if and only if
\begin{align*}
\left(\begin{array}{cc}
        \widetilde{a_1} & \widetilde{b_1} \\
        \widetilde{a_2} & \widetilde{b_2}
      \end{array}
\right)(p,q)\left(\begin{array}{c}
                          \alpha\circ\phi \\
                          \beta\circ\phi
                        \end{array}
\right)(p,q)+\left(\begin{array}{c}
                          \widetilde{\lambda_p} \\
                          \widetilde{\lambda_q}
                        \end{array}
\right)(p,q)=\left(\begin{array}{c}
                          0 \\
                          0
                        \end{array}
\right)
\end{align*}
holds from Lemma \ref{parameter le} and
\begin{align*}
\left(\begin{array}{cc}
  u_p & v_p \\
  u_q & v_q
\end{array}\right)(p,q)\left(\begin{array}{c}
                          \lambda_u \\
                          \lambda_v
                        \end{array}
\right)(\phi(p,q))=\left(\begin{array}{c}
                          \widetilde{\lambda_p} \\
                          \widetilde{\lambda_q}
                        \end{array}
\right)(p,q).
\end{align*}
Thus, it is equivalent to say that there exits a mapping $(\alpha\circ\phi,\beta\circ\phi,\gamma\circ\phi):V\to S^2$ satisfying the creative condition for the sphere family $\mathcal{S}_{(\widetilde{x}, \widetilde{\lambda})}$.
\end{proof}

Let $(\bm{x},\bm{n},\bm{s}): U\to \mathbb{R}^3\times\Delta$ be a framed surface with the basic invariants $\mathcal{A}$, $\mathcal{F}_1$ and $\mathcal{F}_2$. We consider the
rotation of the frame $\{\bm{n},\bm{s},\bm{t}\}$.
Let $\theta:U\rightarrow\mathbb{R}$ be a smooth function, we denote
\begin{align*}
\left(\begin{array}{c}
                          \bm{s}^\theta(u,v) \\
                          \bm{t}^\theta(u,v)
                        \end{array}
\right)=\left(\begin{array}{cc}
  \cos\theta(u,v) & -\sin\theta(u,v) \\
  \sin\theta(u,v) & \cos\theta(u,v)
\end{array}\right)\left(\begin{array}{c}
                          \bm{s}(u,v) \\
                          \bm{t}(u,v)
                        \end{array}
\right).
\end{align*}
Then, $\{\bm{n},\bm{s}^\theta,\bm{t}^\theta\}$ is a moving frame along $\bm{x}$ and $(\bm{x},\bm{n},\bm{s}^\theta)$ is also a framed surface. We denote the basic invariants of $(\bm{x},\bm{n},\bm{s}^\theta)$ by $(\mathcal{A}^\theta, \mathcal{F}_1^\theta, \mathcal{F}_2^\theta)$. There is the following lemma in \cite{fukunagatakahashi1}.

\begin{lemma}\label{rotation le}(\cite{fukunagatakahashi1})
{\rm
Under the above notations, the relation between $(\mathcal{A}, \mathcal{F}_1, \mathcal{F}_2)$ and $(\mathcal{A}^\theta, \mathcal{F}_1^\theta, \mathcal{F}_2^\theta)$ is as follows.
\begin{align*}
\mathcal{A}^\theta=&\left(\begin{array}{cc}
                           a_1\cos\theta-b_1\sin\theta & a_1\sin\theta+b_1\cos\theta \\
                           a_2\cos\theta-b_2\sin\theta & a_2\sin\theta+b_2\cos\theta
                         \end{array}
\right)=\mathcal{A}\left(\begin{array}{cc}
                           \cos\theta & \sin\theta \\
                           -\sin\theta & \cos\theta
                         \end{array}
\right),\\
\mathcal{F}_1^\theta=&\left(\begin{array}{ccc}
                              0 & e_1\cos\theta-f_1\sin\theta & e_1\sin\theta+f_1\cos\theta \\
                              -e_1\cos\theta+f_1\sin\theta & 0 & g_1-\theta_u \\
                              -e_1\sin\theta-f_1\cos\theta  & -g_1+\theta_u & 0
                            \end{array}
\right),\\
\mathcal{F}_2^\theta=&\left(\begin{array}{ccc}
                              0 & e_2\cos\theta-f_2\sin\theta & e_2\sin\theta+f_2\cos\theta \\
                              -e_2\cos\theta+f_2\sin\theta & 0 & g_2-\theta_v \\
                              -e_2\sin\theta-f_2\cos\theta  & -g_2+\theta_v & 0
                            \end{array}
\right).
\end{align*}
}
\end{lemma}

\begin{proposition}
Suppose the sphere family $\mathcal{S}_{(x, \lambda)}$ is creative, then the mapping $(\alpha,\beta,\gamma):U\to\mathbb{R}^3$ satisfies the creative condition related to the frame $\{\bm{n},\bm{s},\bm{t}\}$ if and only if the mapping $(\alpha^\theta,\beta^\theta,\gamma)=(\alpha\cos\theta-\beta\sin\theta,\alpha\sin\theta+\beta\cos\theta,\gamma):U\to\mathbb{R}^3$ satisfies the creative condition related to the frame $\{\bm{n},\bm{s}^\theta,\bm{t}^\theta\}$.
\end{proposition}
\begin{proof}
By Definition \ref{creative}, the mapping $(\alpha,\beta,\gamma):U\to\mathbb{R}^3$ satisfies the creative condition related to the frame $\{\bm{n},\bm{s},\bm{t}\}$ if and only if
\begin{align*}
\mathcal{A}\left(\begin{array}{c}
                   \alpha \\
                   \beta
                 \end{array}
\right)+\left(\begin{array}{c}
                   \lambda_u \\
                   \lambda_v
                 \end{array}
\right)=\left(\begin{array}{c}
                   0 \\
                   0
                 \end{array}
\right)
\end{align*}
and $\alpha^2+\beta^2+\gamma^2=1$. By Lemma \ref{rotation le}, we have
\begin{align*}
\mathcal{A}\left(\begin{array}{c}
                   \alpha \\
                   \beta
                 \end{array}
\right)=\mathcal{A}\left(\begin{array}{cc}
                           \cos\theta & \sin\theta \\
                           -\sin\theta & \cos\theta
                         \end{array}
\right)\left(\begin{array}{cc}
                           \cos\theta & -\sin\theta \\
                           \sin\theta & \cos\theta
                         \end{array}
\right)\left(\begin{array}{c}
                   \alpha \\
                   \beta
                 \end{array}
\right)=\mathcal{A}^\theta\left(\begin{array}{c}
                   \alpha^\theta \\
                   \beta^\theta
                 \end{array}
\right).
\end{align*}
It follows \begin{align*}
\mathcal{A}\left(\begin{array}{c}
                   \alpha \\
                   \beta
                 \end{array}
\right)+\left(\begin{array}{c}
                   \lambda_u \\
                   \lambda_v
                 \end{array}
\right)=\left(\begin{array}{c}
                   0 \\
                   0
                 \end{array}
\right)
\end{align*}
if and only if
\begin{align*}
\mathcal{A}^\theta\left(\begin{array}{c}
                   \alpha^\theta \\
                   \beta^\theta
                 \end{array}
\right)+\left(\begin{array}{c}
                   \lambda_u \\
                   \lambda_v
                 \end{array}
\right)=\left(\begin{array}{c}
                   0 \\
                   0
                 \end{array}
\right).
\end{align*}
Moreover, we can see that $\alpha^2+\beta^2+\gamma^2=1$ is equivalent to $(\alpha^\theta)^2+(\beta^\theta)^2+\gamma^2=1$. Thus, the conclusion holds.
\end{proof}

By Definition \ref{creative}, any family of concentric spheres with smoothly expanding radii is not creative, and it is obvious that such a sphere family does not create an envelope. Moreover, for a sphere family $\mathcal{S}_{(x, \lambda)}$, if the creative condition is satisfied for any $(u,v)\in U$, then $U$ consists of the following five sets.
\begin{align*}
\Sigma_1=&\{(u,v)\in U~|~J_a^2(u,v)+J_b^2(u,v)<J_F^2(u,v)\},\\
\Sigma_2=&\{(u,v)\in U~|~J_a^2(u,v)+J_b^2(u,v)=J_F^2(u,v)\neq0\},\\
\Sigma_3=&\{(u,v)\in U~|~\mathrm{rank}\mathcal{A}=1,~ (a_1^2(u,v)+b_1^2(u,v),a_2^2(u,v)+b_2^2(u,v))=(\lambda_u^2(u,v),\lambda_v^2(u,v))\},\\
\Sigma_4=&\{(u,v)\in U~|~\mathrm{rank}\mathcal{A}=1,~ a_1^2(u,v)+b_1^2(u,v)>\lambda_u^2(u,v)~ \mathrm{or}~a_2^2(u,v)+b_2^2(u,v)>\lambda_v^2(u,v)\},\\
\Sigma_5=&\{(u,v)\in U~|~\mathrm{rank}\mathcal{A}=0,~ \lambda_u(u,v)=\lambda_v(u,v)=0\},
\end{align*}
where
\begin{align*}
J_F(u,v)=\mathrm{det}\mathcal{A}(u,v),~~
J_a(u,v)&=\mathrm{det}\left(
         \begin{array}{cc}
a_1(u,v)& \lambda_u(u,v)\\
      a_2(u,v) &   \lambda_v(u,v) \\
         \end{array}
       \right),~~
J_b(u,v)=\mathrm{det}\left(
         \begin{array}{cc}
           b_1(u,v)&\lambda_u(u,v) \\
           b_2(u,v)&\lambda_v(u,v) \\
         \end{array}
       \right).
\end{align*}

Under the above preparation, Problem \ref{problem} is solved as follows.

\begin{theorem}\label{theorem1}
Let $(\bm{x},\bm{n},\bm{s}): U\to\mathbb{R}^3\times\Delta$
be a framed surface and let
$\lambda: U\to \mathbb{R}_+$ be a positive function.
Then, the following three hold.
\begin{enumerate}
\item[(1)] The sphere family $\mathcal{S}_{(x, \lambda)}$
creates an envelope if and only if $\mathcal{S}_{(x, \lambda)}$ is creative.
\item[(2)] Suppose the sphere family
$\mathcal{S}_{(x, \lambda)}$ creates
an envelope $\bm{f}: U\to \mathbb{R}^3$.     Then, the created envelope
$\bm{f}$ is represented as follows.
\[
\bm{f}(u,v) = \bm{x}(u,v) + \lambda(u,v)\bm{\nu}(u,v),
\]
where $\bm{\nu}: U\to S^2$ is the mapping
defined in Definition \ref{creative}.
\item[(3)] Suppose that the sphere family
$\mathcal{S}_{(x, \lambda)}$ creates
an envelope.
Then, the number of envelopes created by $\mathcal{S}_{(x, \lambda)}$ is
characterized as follows.
\smallskip
\begin{enumerate}
\item[(3-i)] The sphere family
$\mathcal{S}_{(x, \lambda)}$
creates a unique envelope
if the set $\Sigma_2$ is dense in $U$ or the set $\Sigma_3$ is dense in $U$.
\item[(3-ii)] There are exactly two distinct envelopes created by
$\mathcal{S}_{(x, \lambda)}$ if the set $\Sigma_1$ is dense in $U$.
\item[(3-$\infty$)] There are uncountably many distinct
envelopes created by
$\mathcal{S}_{(x, \lambda)}$ if the set $\Sigma_1\cup\Sigma_2\cup\Sigma_3$ is not dense in $U$.
\end{enumerate}
\end{enumerate}
\end{theorem}
\noindent
By the assertion (2) of Theorem \ref{theorem1}, we can reasonably call $\bm{\nu}$ the {\it creator} for an envelope $\bm{f}$ created by $\mathcal{S}_{(x, \lambda)}$. Actually, by Definition \ref{envelope}, $(\bm{f},\bm{\nu}):U\rightarrow\mathbb{R}^3\times S^2$ is a Legendre surface. Moreover, since $\bm{\nu}(u,v)=\alpha(u,v)\bm{s}(u,v)+\beta(u,v)\bm{t}(u,v)+\gamma(u,v)\bm{n}(u,v)$, we denote
\begin{enumerate}
\item[(i)]  $\bm{\omega}(u,v)=\bm{n}(u,v)$ if $\gamma(u,v)\equiv0$,
\item[(ii)] $\bm{\omega}(u,v)=\frac{1}{\sqrt{\gamma^2(u,v)+\alpha^2(u,v)}}
\left(-\gamma(u,v)\bm{s}(u,v)+\alpha(u,v)\bm{n}(u,v)\right)$ if $\gamma(u,v)\neq0$.
\end{enumerate}
Then $(\bm{f},\bm{\nu},\bm{\omega}):U\rightarrow\mathbb{R}^3\times\Delta$ is a framed surface. We have the following proposition by a direct calculation.
\begin{proposition}\label{finvariants}
Let $(\bm{f},\bm{\nu},\bm{\omega}):U\rightarrow\mathbb{R}^3\times\Delta$ be a framed surface and $\bm{\mu}(u,v)=\bm{\nu}(u,v)\times \bm{\omega}(u,v)$. Then $\{\bm{\nu}(u,v),\bm{\omega}(u,v),\bm{\mu}(u,v)\}$ is a moving frame along $\bm{f}(u,v)$, and we have the following systems of differential equations:
\begin{align*}
\left(
  \begin{array}{c}
    \bm{f}_u \\
    \bm{f}_v \\
  \end{array}
\right)&=\left(
          \begin{array}{cc}
            a_{f1} & b_{f1} \\
            a_{f2} & b_{f2} \\
          \end{array}
        \right)
        \left(
          \begin{array}{c}
            \bm{\omega} \\
            \bm{\mu} \\
          \end{array}
        \right),
\\
\left(
  \begin{array}{c}
    \bm{\nu}_u \\
    \bm{\omega}_u \\
    \bm{\mu}_u \\
  \end{array}
\right)&=\left(
          \begin{array}{ccc}
            0 & e_{f1} & f_{f1} \\
            -e_{f1} & 0 & g_{f1} \\
            -f_{f1} & -g_{f1} & 0  \\
          \end{array}
        \right)
        \left(
          \begin{array}{c}
            \bm{\nu} \\
            \bm{\omega} \\
            \bm{\mu} \\
          \end{array}
        \right),\\
\left(
  \begin{array}{c}
    \bm{\nu}_v \\
    \bm{\omega}_v \\
    \bm{\mu}_v \\
  \end{array}
\right)&=\left(
          \begin{array}{ccc}
            0 & e_{f2} & f_{f2} \\
            -e_{f2} & 0 & g_{f2} \\
            -f_{f2} & -g_{f2} & 0  \\
          \end{array}
        \right)
        \left(
          \begin{array}{c}
            \bm{\nu} \\
            \bm{\omega} \\
            \bm{\mu} \\
          \end{array}
        \right),
\end{align*}
where $a_{fi}$, $b_{fi}$, $e_{fi}$, $f_{fi}$, $g_{fi}$ (i=1,2) are the basic invariants of the framed surface $(\bm{f},\bm{\nu},\bm{\omega})$.
\end{proposition}

\section{Proof of Theorem \ref{theorem1}}\label{section2}
\subsection{Proof of the assertion (1) of Theorem \ref{theorem1}}
\label{subsection2.1}
\begin{proof}
Suppose that a sphere family $\mathcal{S}_{(x, \lambda)}$ is creative. By definition, there exists a smooth mapping $\bm{\nu}:U\to S^2$ such that the following identities hold for any $(u,v)\in U$.
\begin{align*}
 \bm{x}_u(u,v)\cdot\bm{\nu}(u,v)+\lambda_u(u,v)=0,~~~~~~
 \bm{x}_v(u,v)\cdot\bm{\nu}(u,v)+\lambda_v(u,v)=0.
\end{align*}
Set
\begin{align*}
\bm{f}(u,v)=\bm{x}(u,v)+\lambda(u,v)\bm{\nu}(u,v),
\end{align*}
it follows $\bm{f}(u,v)\in\mathcal{S}_{(x, \lambda)}$ since
\begin{align*}
\left(\bm{f}(u,v)-\bm{x}(u,v)\right)\cdot\left(\bm{f}(u,v)-\bm{x}(u,v)\right)=\lambda^2(u,v).
\end{align*}
Moreover, since
\begin{align*}
\bm{f}_u(u,v)&=\bm{x}_u(u,v)+\lambda_u(u,v)\bm{\nu}(u,v)+\lambda(u,v)\bm{\nu}_u(u,v),\\
\bm{f}_v(u,v)&=\bm{x}_v(u,v)+\lambda_v(u,v)\bm{\nu}(u,v)+\lambda(u,v)\bm{\nu}_v(u,v),
\end{align*}
we have
\begin{align*}
&\bm{f}_u(u,v)\cdot\left(\bm{f}(u,v)-\bm{x}(u,v)\right)\\
&=\left(\bm{x}_u(u,v)+\lambda_u(u,v)\bm{\nu}(u,v)+\lambda(u,v)\bm{\nu}_u(u,v)\right)\cdot\lambda(u,v)\bm{\nu}(u,v)\\
&=\lambda(u,v)\bm{x}_u(u,v)\cdot\bm{\nu}(u,v)+\lambda(u,v)\lambda_u(u,v)\\
&=0
\end{align*}
and
\begin{align*}
&\bm{f}_v(u,v)\cdot\left(\bm{f}(u,v)-\bm{x}(u,v)\right)\\
&=\left(\bm{x}_v(u,v)+\lambda_v(u,v)\bm{\nu}(u,v)+\lambda(u,v)\bm{\nu}_v(u,v)\right)\cdot\lambda(u,v)\bm{\nu}(u,v)\\
&=\lambda(u,v)\bm{x}_v(u,v)\cdot\bm{\nu}(u,v)+\lambda(u,v)\lambda_v(u,v)\\
&=0.
\end{align*}
Thus, $\bm{f}$ is an envelope created by $\mathcal{S}_{(x, \lambda)}$.

Conversely, suppose that the sphere family $\mathcal{S}_{(x, \lambda)}$ creates an envelope $\bm{f}:U\to \mathbb{R}^3$. By definition, we have $\bm{f}\in\mathcal{S}_{(x, \lambda)}$ and
\begin{align*}
\bm{f}_u(u,v)\cdot\left(\bm{f}(u,v)-\bm{x}(u,v)\right)=\bm{f}_v(u,v)\cdot\left(\bm{f}(u,v)-\bm{x}(u,v)\right)=0.
\end{align*}
The condition $\bm{f}\in\mathcal{S}_{(x, \lambda)}$ implies that there exists a smooth mapping $\bm{\nu}:U\to S^2$ such that
\begin{align*}
\bm{f}(u,v)=\bm{x}(u,v)+\lambda(u,v)\bm{\nu}(u,v).
\end{align*}
Then, since
\begin{align*}
\bm{f}_u(u,v)&=\bm{x}_u(u,v)+\lambda_u(u,v)\bm{\nu}(u,v)+\lambda(u,v)\bm{\nu}_u(u,v),\\
\bm{f}_v(u,v)&=\bm{x}_v(u,v)+\lambda_v(u,v)\bm{\nu}(u,v)+\lambda(u,v)\bm{\nu}_v(u,v),
\end{align*}
we have
\begin{align*}
0&=\bm{f}_u(u,v)\cdot\left(\bm{f}(u,v)-\bm{x}(u,v)\right)\\
&=\left(\bm{x}_u(u,v)+\lambda_u(u,v)\bm{\nu}(u,v)+\lambda(u,v)\bm{\nu}_u(u,v)\right)\cdot\lambda(u,v)\bm{\nu}(u,v)\\
&=\lambda(u,v)\left(\bm{x}_u(u,v)\cdot\bm{\nu}(u,v)+\lambda_u(u,v)\right)
\end{align*}
and
\begin{align*}
0&=\bm{f}_v(u,v)\cdot\left(\bm{f}(u,v)-\bm{x}(u,v)\right)\\
&=\left(\bm{x}_v(u,v)+\lambda_v(u,v)\bm{\nu}(u,v)+\lambda(u,v)\bm{\nu}_v(u,v)\right)\cdot\lambda(u,v)\bm{\nu}(u,v)\\
&=\lambda(u,v)\left(\bm{x}_v(u,v)\cdot\bm{\nu}(u,v)+\lambda_v(u,v)\right).
\end{align*}
Since $\lambda(u,v)$ is positive for any $(u,v)\in U$, it follows
\begin{align*}
 \bm{x}_u(u,v)\cdot\bm{\nu}(u,v)+\lambda_u(u,v)=0,~~~~~~
 \bm{x}_v(u,v)\cdot\bm{\nu}(u,v)+\lambda_v(u,v)=0.
\end{align*}
Thus, the sphere family $\mathcal{S}_{(x, \lambda)}$ is creative.
\end{proof}
\subsection{Proof of the assertion (2) of Theorem \ref{theorem1}}
\label{subsection2.2}
\begin{proof}
The proof of the assertion (1) given in Subsection \ref{subsection2.1} proves
the assertion (2) as well.
\end{proof}

\subsection{Proof of the assertion (3) of Theorem \ref{theorem1}}
\label{subsection2.3}
\subsubsection{Proof of (3-i)}
\begin{proof}
Since the sphere family $\mathcal{S}_{(x, \lambda)}$ is creative, then there exist three functions $\alpha(u,v)$, $\beta(u,v)$ and $\gamma(u,v)$ satisfying
\begin{align*}
&a_1(u,v)\alpha(u,v)+b_1(u,v)\beta(u,v)+\lambda_u(u,v)=0,\\
&a_2(u,v)\alpha(u,v)+b_2(u,v)\beta(u,v)+\lambda_v(u,v)=0,\\
&\alpha^2(u,v)+\beta^2(u,v)+\gamma^2(u,v)=1.
\end{align*}
Suppose that the set $\Sigma_2$ is dense in $U$ or the set $\Sigma_3$ is dense in $U$. Then
under considering continuity of the functions $(u,v)\mapsto\alpha(u,v)$, $(u,v)\mapsto\beta(u,v)$ and $(u,v)\mapsto\gamma(u,v)$, it follows from the assumption that $\alpha^2(u,v)+\beta^2(u,v)=1$ hold for any $(u,v)\in U$. Meanwhile, $\alpha(u,v)$, $\beta(u,v)$ are uniquely determined for any $(u,v)\in U$.
By Definition \ref{creative}, we have
\begin{align*}
\bm{\nu}(u,v)=\alpha(u,v)\bm{s}(u,v)+\beta(u,v)\bm{t}(u,v)+\gamma(u,v)\bm{n}(u,v)=\alpha(u,v)\bm{s}(u,v)+\beta(u,v)\bm{t}(u,v).
\end{align*}
Therefore, the creator $\bm{\nu}(u,v)$ is uniquely determined for any $(u,v)\in U$. By the assertion (2) of Theorem \ref{theorem1}, the created envelope $\bm{f}(u,v)=\bm{x}(u,v)+\lambda(u,v)\bm{\nu}(u,v)$ must be unique.
\end{proof}

\subsubsection{Proof of (3-ii)}
\begin{proof}
Let $(u_0,v_0)$ be a point of $\Sigma_1$. Then, $J_a^2(u_0,v_0)+J_b^2(u_0,v_0)<J_F^2(u_0,v_0)$ holds. It follows that there must exist an open set $\widetilde{U}\subseteq U$ such that $(u_0,v_0)\in\widetilde{U}$ and $\alpha^2(u,v)+\beta^2(u,v)<1$ for any $(u,v)\in\widetilde{U}$. Note that $\alpha^2(u,v)+\beta^2(u,v)+\gamma^2(u,v)=1,$ $\gamma(u,v)$ can take two distinct values. Namely,
\begin{align*}
\gamma_+(u,v)&=\sqrt{1-\alpha^2(u,v)-\beta^2(u,v)},\\ \gamma_-(u,v)&=-\sqrt{1-\alpha^2(u,v)-\beta^2(u,v)}.
\end{align*}
Hence, for any $(u,v)\in\widetilde{U}$, there exist exactly two distinct unit
vectors $\bm{\nu}_+(u,v)$, $\bm{\nu}_-(u,v)$ corresponding to $\gamma_+(u,v)$ and $\gamma_-(u,v)$, respectively. Since the set $\Sigma_1$ is dense in $U$, it follows that there exist exactly two distinct creators $\bm{\nu}_+:U\rightarrow S^2$ and $\bm{\nu}_-:U\rightarrow S^2$ for the sphere family $\mathcal{S}_{(x, \lambda)}$. Thus, by the assertion (2) of Theorem \ref{theorem1}, the sphere family $\mathcal{S}_{(x, \lambda)}$ must create two distinct envelopes.
\end{proof}
\subsubsection{Proof of (3-$\infty$)}
\begin{proof}
Suppose that the set $\Sigma_1\cup\Sigma_2\cup\Sigma_3$ is not dense in $U$. It follows that there exits an open set $\widetilde{U}\subseteq\Sigma_4\cup\Sigma_5\subseteq U$. For any $(u,v)\in\widetilde{U}$,
solve the following equalities
\begin{align*}
&a_1(u,v)\alpha(u,v)+b_1(u,v)\beta(u,v)+\lambda_u(u,v)=0,\\
&a_2(u,v)\alpha(u,v)+b_2(u,v)\beta(u,v)+\lambda_v(u,v)=0,
\end{align*}
the solution $(\alpha(u,v),\beta(u,v))$ forms a line segment in the unit disk if $(u,v)\in\Sigma_4$, or the entire unit disk if $(u,v)\in\Sigma_5$. Then $\bm{\nu}(u,v)$ can take many distinct unit vectors for any $(u,v)\in\widetilde{U}$.
Take one point $(u_0,v_0)\in\widetilde{U}$ and denote the $\bm{\nu}$ for the envelope $\bm{f}_0$ by $\bm{\nu}_0$.
Then, by using the standard technique on
bump functions, we may construct uncountably many distinct creators $\bm{\nu}_\epsilon:U\to S^2$ $(\epsilon\in A)$ such that the following conditions $(a),(b)$ and $(c)$ hold, where $A$ is a set consisting uncountably many elements such that $0\notin A$.
\begin{enumerate}
\item[(a)] For any $(u,v)\in U-\widetilde{U}$ and any $\epsilon\in A$,
the equality
$\bm{\nu}_\epsilon(u,v)=\bm{\nu}_0(u,v)$ holds.
\item[(b)] For any $\epsilon\in A$, the property
$\bm{\nu}_\epsilon(u_0,v_0)\ne \bm{\nu}_0(u_0,v_0)$ holds.
\item[(c)] For any two distinct $\epsilon_1, \epsilon_2\in A$, the property
$\bm{\nu}_{\epsilon_1}(u_0,v_0)\ne
\bm{\nu}_{\epsilon_2}(u_0,v_0)$ holds.
\end{enumerate}
Therefore, by the assertion (2) of Theorem \ref{theorem1}, there are uncountably many distinct envelopes $\bm{f}_\epsilon$ created by $\mathcal{S}_{(x, \lambda)}$.
\end{proof}

\section{Examples}\label{section3}
\begin{example}\label{example3}
We examine Example \ref{example1} by applying Theorem \ref{theorem1}.
In Example \ref{example1}, $(\bm{x},\bm{n},\bm{s}): \mathbb{R}^2\to \mathbb{R}^3\times\Delta$
is a framed surface given by $\bm{x}(u,v)=\left(u^2, u^3, v\right)$, $\bm{n}(u,v)=(1/\sqrt{9u^2+4})(-3u,2,0)$ and $\bm{s}=(0,0,1)$.
Moreover, the radius function
$\lambda: \mathbb{R}^2\to \mathbb{R}_+$ is defined by $\lambda(u,v)=1$.
Thus,
\[
\lambda_{u}(u,v)=0,~~~~~~\lambda_{v}(u,v)=0.
\]
Therefore, the unit vector $\bm{\nu}(u,v)\in S^2$ satisfying
\begin{align*}
\bm{x}_u(u,v)\cdot\bm{\nu}(u,v)+\lambda_u(u,v)=0,~~~~~~
\bm{x}_v(u,v)\cdot\bm{\nu}(u,v)+\lambda_v(u,v)=0
\end{align*}
exists and it must have the form
\[
\bm{\nu}(u,v)=\pm \bm{n}(u,v)=\pm \frac{1}{\sqrt{9u^2+4}}(-3u,2,0).
\]
Hence, by the assertion (1) of Theorem \ref{theorem1}, the sphere family
$\mathcal{S}_{(x, \lambda)}$ creates an envelope
$\bm{f}: \mathbb{R}^2\to \mathbb{R}^3$.
By the assertion (2) of Theorem \ref{theorem1}, $\bm{f}$ is parametrized as follows.
\begin{align*}
\bm{f}(u,v)  = & ~\bm{x}(u,v)+\lambda(u,v)\bm{\nu}(u,v) \\
 = & \left(u^2, u^3, v\right) \pm (1/\sqrt{9u^2+4})\left(-3u,2,0\right) \\
 = &
\left(
u^2\pm \frac{3u}{\sqrt{9u^2+4}}, u^3\mp \frac{2}{\sqrt{9u^2+4}}, v
\right).
\end{align*}
Moreover, by calculation, we have
\begin{align*}
\mathcal{A}&=\left(
          \begin{array}{cc}
            a_1 & b_1 \\
            a_2 & b_2 \\
          \end{array}
        \right)
        =\left(
          \begin{array}{cc}
            0 & u\sqrt{9u^2+4}\\
            1 & 0 \\
          \end{array}
        \right)
\end{align*}
and $J_a=J_b=0$. It follows the set $\Sigma_1$ is dense in $\mathbb{R}^2.$
Finally, by the assertion (3-ii) of Theorem \ref{theorem1},
the number of distinct envelopes created by the
sphere family $\mathcal{S}_{(x, \lambda)}$ is exactly two.
Therefore, Theorem \ref{theorem1} reveals that
the set $\mathcal{D}$ calculated in Example \ref{example1} is
certainly the union of the unit cylindrical surface and the set of two envelopes of
$\mathcal{S}_{(x, \lambda)}$.
\end{example}

\begin{example}\label{example41}
We examine (1) of Example \ref{example2} by applying Theorem \ref{theorem1}.
In (1) of Example \ref{example2}, $\bm{x}: \mathbb{R}^2
{\color{black}\setminus\{(0,0)\}}
\to \mathbb{R}^3$
is given by $\bm{x}(u,v)=\left(u, v, 0\right)$.
If we define two unit vectors $\bm{n}(u,v)=(0, 0, 1)$ and $\bm{s}(u,v)=(1,0,0)$, then
$(\bm{x},\bm{n},\bm{s}): \mathbb{R}^2
{\color{black}\setminus\{(0,0)\}}
\to \mathbb{R}^3\times\Delta$ is a framed surface.
By definition, we have
\begin{align*}
\mathcal{A}&=\left(
          \begin{array}{cc}
            a_1 & b_1 \\
            a_2 & b_2 \\
          \end{array}
        \right)
        =\left(
          \begin{array}{cc}
            1 & 0\\
            0 & 1 \\
          \end{array}
        \right).
\end{align*}
On the other hand, the radius function
$\lambda: \mathbb{R}^2
{\color{black}\setminus\{(0,0)\}}
\to \mathbb{R}_+$ is defined by
$\lambda(u,v)=\sqrt{u^2+v^2}$ in this example.
Thus, the creative condition
\begin{align*}
\bm{x}_u(u,v)\cdot\bm{\nu}(u,v)+\lambda_u(u,v)=0,~~~~~~
\bm{x}_v(u,v)\cdot\bm{\nu}(u,v)+\lambda_v(u,v)=0
\end{align*}
simply becomes
\begin{align*}
\bm{x}_u(u,v)\cdot\bm{\nu}(u,v)+\frac{u}{\sqrt{u^2+v^2}}=0,~~~~~~
\bm{x}_v(u,v)\cdot\bm{\nu}(u,v)+\frac{v}{\sqrt{u^2+v^2}}=0
\end{align*}
in this case.
If we set $\bm{\nu}(u,v)=(-u/\sqrt{u^2+v^2},-v/\sqrt{u^2+v^2},0)$, then the above equality holds for any
$(u,v)\in\mathbb{R}^2
{\color{black}\setminus\{(0,0)\}}
$.
Thus, by the assertion (1) of Theorem \ref{theorem1}, the sphere
family $\mathcal{S}_{(x, \lambda)}$ creates an envelope.
By the assertion (2) of Theorem \ref{theorem1},
the parametrization of the created envelope is
\begin{align*}
\bm{f}(u,v)=&~\bm{x}(u,v)+\lambda(u,v)\bm{\nu}(u,v)\\
=& \left(u, v, 0\right)+\sqrt{u^2+v^2}\left(-\frac{u}{\sqrt{u^2+v^2}},-\frac{v}{\sqrt{u^2+v^2}},0\right)\\ =&\left(0, 0, 0\right).
\end{align*}
Finally, notice that for any
$(u,v)\in\mathbb{R}^2
{\color{black}\setminus\{(0,0)\}}
$
the equality $J_a^2(u,v)+J_b^2(u,v)=J_F^2(u,v)=1$ always
holds. This means the set $\Sigma_2$ is dense in $\mathbb{R}^2
{\color{black}\setminus\{(0,0)\}}
$.
Thus, by the assertion (3-i) of Theorem \ref{theorem1},
the origin $(0,0,0)$
{\color{black}of $\mathbb{R}^3$}
is the unique envelope created by
$\mathcal{S}_{(x, \lambda)}$. The sphere family $\mathcal{S}_{(x, \lambda)}$
and the candidates of its envelope are depicted in
Figure \ref{figure_example5}.
\begin{figure}[h]
\begin{center}
\includegraphics[width=7cm]
{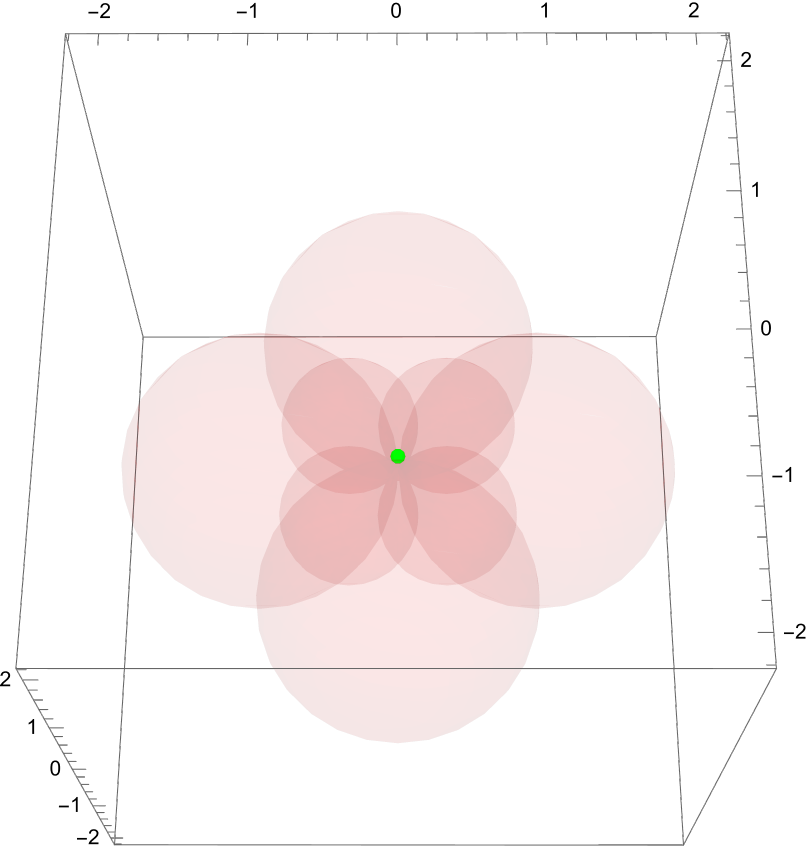}
\caption{The sphere family $\mathcal{S}_{(x, \lambda)}$
and the candidates of its envelope.
}
\label{figure_example5}
\end{center}
\end{figure}
\end{example}

\begin{example}\label{example5}
Theorem \ref{theorem1} can be applied also
to (2) of Example \ref{example2} as follows.
In this example,
$\bm{x}(u,v)=(u,v,0)$ and $\lambda(u,v)=1$, where $(u,v)\in\mathbb{R}^2$.
Thus, we may also set $\bm{n}(u,v)=(0, 0, 1)$ and $\bm{s}(u,v)=(1,0,0)$.
It follows
\begin{align*}
\mathcal{A}&=\left(
          \begin{array}{cc}
            a_1 & b_1 \\
            a_2 & b_2 \\
          \end{array}
        \right)
        =\left(
          \begin{array}{cc}
            1 & 0\\
            0 & 1 \\
          \end{array}
        \right).
\end{align*}
Since the radius function $\lambda$ is a constant function,
the creative condition
\begin{align*}
\bm{x}_u(u,v)\cdot\bm{\nu}(u,v)+\lambda_u(u,v)=0,~~~~~~
\bm{x}_v(u,v)\cdot\bm{\nu}(u,v)+\lambda_v(u,v)=0
\end{align*}
simply becomes
\begin{align*}
\bm{x}_u(u,v)\cdot\bm{\nu}(u,v)=\bm{x}_v(u,v)\cdot\bm{\nu}(u,v)=0
\end{align*}
in this case.
Thus, for any $(u,v)\in \mathbb{R}^2$,
the creative condition is satisfied if and only if
$\bm{\nu}(u,v)=\pm \bm{n}(u,v)=\pm(0, 0, 1)$.
Hence, by the assertion (1) of Theorem \ref{theorem1}, the sphere
family $\mathcal{S}_{(x, \lambda)}$ creates an envelope.
By the assertion (2) of Theorem \ref{theorem1},
the parametrization of the created envelope is
\[
\bm{f}(u,v)=\bm{x}(u,v)+\lambda(u,v)\bm{\nu}(u,v)
= \left(u,v, 0\right)\pm\left(0,0, 1\right) =\left(u,v, \pm1\right).
\]
Finally, it is easily see that $J_a^2(u,v)+J_b^2(u,v)<J_F^2(u,v)$ always
holds for any $(u,v)\in\mathbb{R}^2$. This means the set $\Sigma_1$ is dense in $\mathbb{R}^2$. By the assertion (3-ii) of Theorem \ref{theorem1},
the number of envelope created by
$\mathcal{S}_{(x, \lambda)}$ is exactly two. The sphere family $\mathcal{S}_{(x, \lambda)}$
and the candidates of its envelope are depicted in
Figure \ref{figure_example6}.
\begin{figure}[h]
\begin{center}
\includegraphics[width=7cm]
{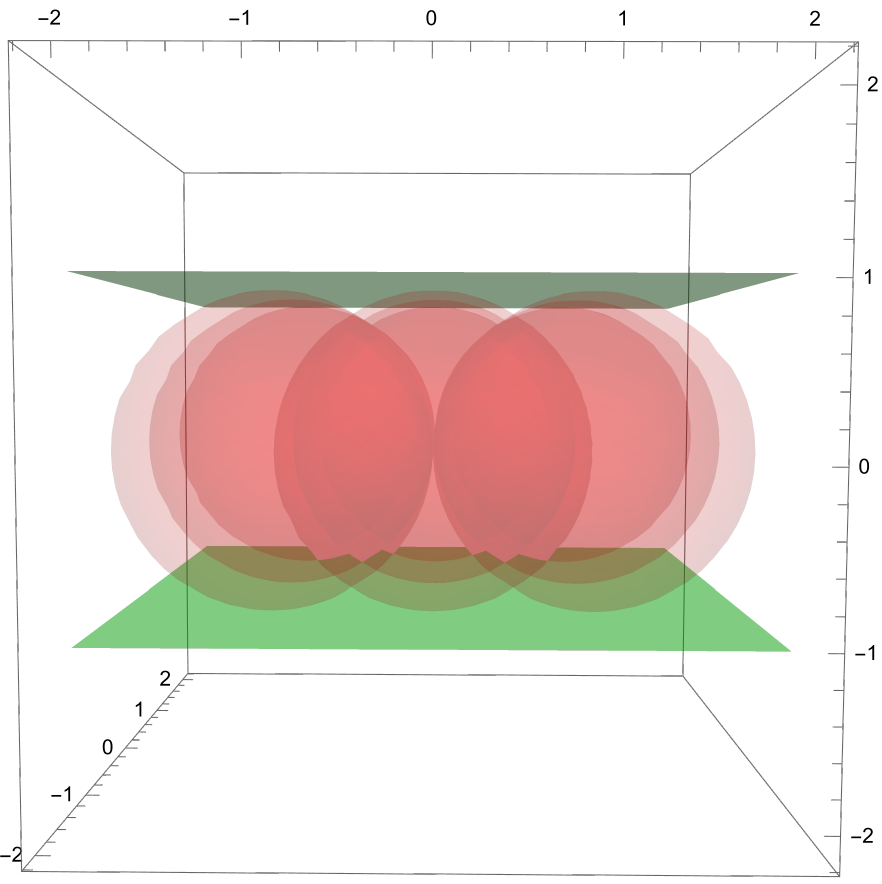}
\caption{The sphere family $\mathcal{S}_{(x, \lambda)}$
and the candidates of its envelope.
}
\label{figure_example6}
\end{center}
\end{figure}
\end{example}

\begin{example}\label{example6}
Theorem \ref{theorem1} can be applied even
to (3) of Example \ref{example2} as follows.
In this example,
$\bm{x}(u,v)=(0,0, 0)$ and $\lambda(u,v)=1$.
Thus, we can take every mapping $(\bm{n},\bm{s}):\mathbb{R}^2\to \Delta$ such that $(\bm{x},\bm{n},\bm{s}): \mathbb{R}^2\to \mathbb{R}^3\times\Delta$ is a framed surface.  In particular,
since the radius function $\lambda$ is a constant function,
the creative condition
\begin{align*}
\bm{x}_u(u,v)\cdot\bm{\nu}(u,v)+\lambda_u(u,v)=0,~~~~~~
\bm{x}_v(u,v)\cdot\bm{\nu}(u,v)+\lambda_v(u,v)=0
\end{align*}
simply becomes
\[
0 = 0,~~0 = 0
\]
in this case.
Thus, for any $\bm{\nu}: \mathbb{R}^2\to S^2$,
the creative condition is satisfied.
Hence, by the assertion (1) of Theorem \ref{theorem1}, the sphere
family $\mathcal{S}_{(x, \lambda)}$ creates an envelope.
By the assertion (2) of Theorem \ref{theorem1},
the parametrization of the created envelope is
\[
\bm{f}(u,v)=\bm{x}(u,v)+\lambda(u,v)\bm{\nu}(u,v)
= \left(0,0, 0\right)+\bm{\nu}(u,v)=\bm{\nu}(u,v).
\]
Finally, since $J_F(u,v)=0$ and $(a_1^2+b_1^2, a_2^2+b_2^2)=(0,0)$ are always hold for any $(u,v)\in\mathbb{R}^2$. This means the set $\Sigma_5$ is dense in $\mathbb{R}^2$. By the assertion (3-$\infty$) of Theorem \ref{theorem1},
there are uncountably many distinct envelopes created by
$\mathcal{S}_{(x, \lambda)}$.
\end{example}

\begin{example}\label{example10}
Let $\bm{x}: \mathbb{R}^2\setminus\mathbb{L}\to \mathbb{R}^3$ be the mapping defined by
$\bm{x}(u,v)=(0, 0, v)$, where $\mathbb{L}=\{(u,v)\in\mathbb{R}^2|v=0\}$. If we take $\bm{n}(u,v)=(1, 0, 0)$ and $\bm{s}(u,v)=(0,1,0)$, then
$(\bm{x},\bm{n},\bm{s}): \mathbb{R}^2\setminus\mathbb{L}\to \mathbb{R}^3\times\Delta$ is a framed surface. Let $\lambda: \mathbb{R}^2\setminus\mathbb{L}\to \mathbb{R}_+$ be the function defined by $\lambda(u,v)=v$.
We calculate that
\begin{align*}
\mathcal{A}&=\left(
          \begin{array}{cc}
            a_1 & b_1 \\
            a_2 & b_2 \\
          \end{array}
        \right)
        =\left(
          \begin{array}{cc}
            0 & 0\\
            0 & 1 \\
          \end{array}
        \right).
\end{align*}
Since $\lambda_u(u,v)=0$ and $\lambda_v(u,v)=1$, it follows $\bm{\nu}(u,v)=(0,0,-1)$ from the creative condition. Thus, by the assertion (2) of Theorem \ref{theorem1}, an envelope $\bm{f}$
created by $\mathcal{S}_{(x, \lambda)}$ is
parametrized as follows.
\begin{align*}
\bm{f}(u,v)  = & ~\bm{x}(u,v)+\lambda(u,v)\bm{\nu}(u,v) \\
 = &
(0, 0, v)+v(0,0,-1)  \\
 = &
(0,0, 0).
\end{align*}
Finally, since $J_F(u,v)=0$, $a_1^2(u,v)+b_1^2(u,v)=\lambda_u^2(u,v)=0$ and $a_2^2(u,v)+b_2^2(u,v)=\lambda_v^2(u,v)=1$ always hold for any $(u,v)\in \mathbb{R}^2$. This means the set $\Sigma_3$ is dense in $\mathbb{R}^2\setminus\mathbb{L}$. By the assertion (3-i) of Theorem \ref{theorem1}, it follows that the
sphere family $\mathcal{S}_{(x, \lambda)}$ creates a unique envelope $(0,0,0)$. The sphere family $\mathcal{S}_{(x, \lambda)}$
and the candidates of its envelope are depicted in
Figure \ref{figure_example8}.
\begin{figure}[h]
\begin{center}
\includegraphics[width=7cm]
{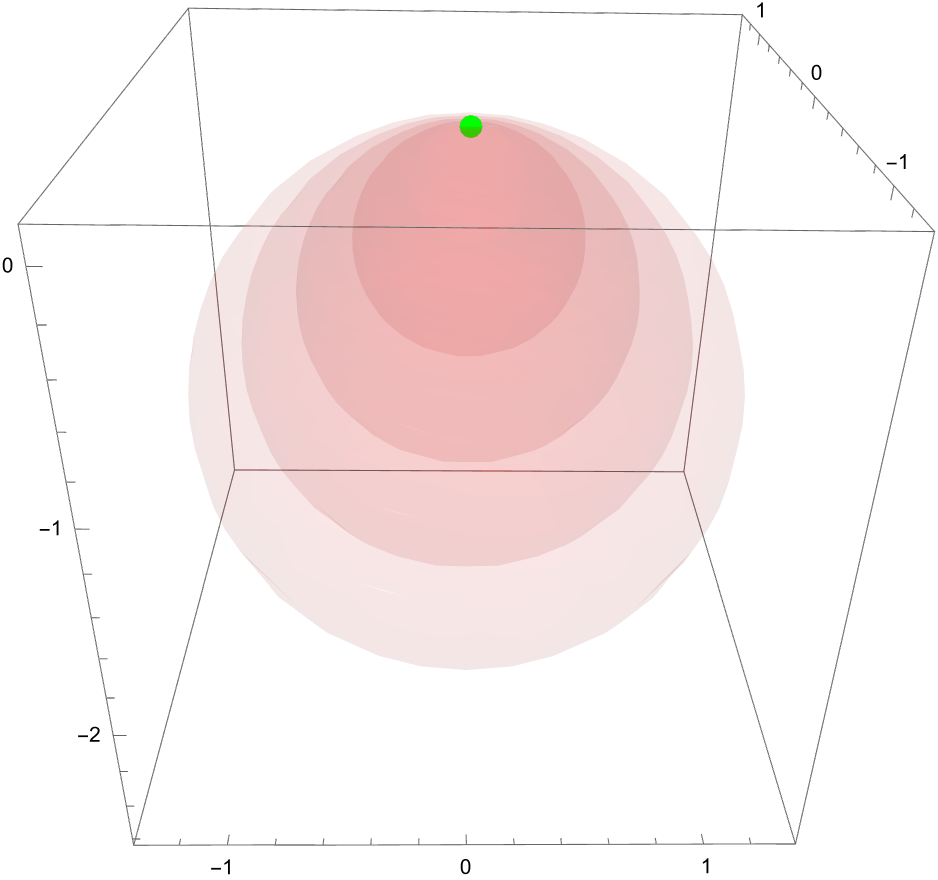}
\caption{The sphere family $\mathcal{S}_{(x, \lambda)}$
and the candidates of its envelope.
}
\label{figure_example8}
\end{center}
\end{figure}
\end{example}

\begin{example}
This example is almost the same as Example \ref{example10}. The difference from Example \ref{example10} is only the radius function $\lambda(u,v)$. Here $\lambda(u,v)=v/2$. Then, from calculations in Example \ref{example10}, we have
\begin{align*}
\mathcal{A}&=\left(
          \begin{array}{cc}
            a_1 & b_1 \\
            a_2 & b_2 \\
          \end{array}
        \right)
        =\left(
          \begin{array}{cc}
            0 & 0\\
            0 & 1 \\
          \end{array}
        \right).
\end{align*}
Since $\lambda_u(u,v)=0$ and $\lambda_v(u,v)=1/2$, by the creative condition
\begin{align*}
\bm{x}_u(u,v)\cdot\bm{\nu}(u,v)+\lambda_u(u,v)=0,~~~~~~
\bm{x}_v(u,v)\cdot\bm{\nu}(u,v)+\lambda_v(u,v)=0,
\end{align*}
we obtain $\bm{\nu}(u,v)\cdot(0,0,1)=-1/2$. Thus, for any function $\theta:\mathbb{R}^2\setminus\mathbb{L}\to \mathbb{R}$, $$\bm{\nu}(u,v)=\bigg(\frac{\sqrt{3}}{2}\cos\theta(u,v),\frac{\sqrt{3}}{2}\sin\theta(u,v),-\frac{1}{2}\bigg)$$ is the creator. By the assertion (2) of Theorem \ref{theorem1}, an envelope $\bm{f}$
created by $\mathcal{S}_{(x, \lambda)}$ is
parametrized as follows.
\begin{align*}
\bm{f}(u,v)  = &~ x(u,v)+\lambda(u,v)\bm{\nu}(u,v) \\
 = &
(0, 0, v)+\frac{1}{2}v\bigg(\frac{\sqrt{3}}{2}\cos\theta(u,v),\frac{\sqrt{3}}{2}\sin\theta(u,v),-\frac{1}{2}\bigg)  \\
 = &
\bigg(\frac{\sqrt{3}}{4}v\cos\theta(u,v),\frac{\sqrt{3}}{4}v\sin\theta(u,v), \frac{3}{4}v\bigg).
\end{align*}
Finally, since $J_F(u,v)=0$ and $a_2^2(u,v)+b_2^2(u,v)>\lambda_v^2(u,v)$ are always hold for any $(u,v)\in\mathbb{R}^2\setminus\mathbb{L}$. This means the set $\Sigma_4$ is dense in $\mathbb{R}^2\setminus\mathbb{L}$. By the assertion (3-$\infty$) of Theorem \ref{theorem1},
there are uncountably many distinct envelopes created by
$\mathcal{S}_{(x, \lambda)}$.
\end{example}

\section{Alternative definitions}\label{section4}
In Definition \ref{envelope} of Section \ref{section1},
the definition of envelope created by the sphere family
is given.
We call the set consisting of the images of envelopes defined
in Definition \ref{envelope}
a \textit{$f$ envelope} (denoted by $\bm{f}$). In \cite{brucegiblin},
 an alternative definition
(called \textit{$\mathcal{D}$ envelope})
is given as follows.
\begin{definition}[$\mathcal{D}$ envelope \cite{brucegiblin}]
{\rm
Let $\bm{x}: U\to \mathbb{R}^3$, $\lambda: U\to \mathbb{R}_+$
be a framed base surface and a positive function respectively.
Set
\[
F(x_1, x_2, x_3, u,v)=||(x_1, x_2, x_3)-\bm{x}(u,v)||^2-\lambda^2(u,v).
\]
Then, the following set is called
the \textit{$\mathcal{D}$ envelope} created by the sphere family
$\mathcal{S}_{(x, \lambda)}$ and is denoted by $\mathcal{D}$.
\[
\left\{(x_1, x_2, x_3)\in \mathbb{R}^3\, |\,
\exists (u,v)\in U \mbox{ s.t. } F(x_1, x_2, x_3, u,v)=\frac{\partial F}{\partial u}(x_1, x_2, x_3, u,v)=\frac{\partial F}{\partial v}(x_1, x_2, x_3, u,v)=0
\right\}.
\]
}
\end{definition}
In Section \ref{section3}, we show a relationship $\bm{f}\subseteq\mathcal{D}$ for Example \ref{example3}.
In this section, we study more precise relationships between $f$ envelope and $\mathcal{D}$ envelope. Suppose the sphere family $\mathcal{S}_{(x, \lambda)}$ is creative. If $\Sigma_4\neq\emptyset$, for a point $(u_0,v_0)\in\Sigma_4$, let $\widetilde{\bm{\nu}}(u_0,v_0)$ be a unit vector represented as
\begin{align*}
\widetilde{\bm{\nu}}(u_0,v_0)=\widetilde{\alpha}(u_0,v_0) \bm{s}(u_0,v_0)+\widetilde{\beta}(u_0,v_0) \bm{t}(u_0,v_0)+\widetilde{\gamma}(u_0,v_0) \bm{n}(u_0,v_0),
\end{align*}
where $\widetilde{\alpha}(u_0,v_0)$, $\widetilde{\beta}(u_0,v_0)$ and $\widetilde{\gamma}(u_0,v_0)$ satisfying
\begin{align*}
&a_1(u_0,v_0)\widetilde{\alpha}(u_0,v_0)+b_1(u_0,v_0)\widetilde{\beta}(u_0,v_0)+\lambda_u(u_0,v_0)=0,\\
&a_2(u_0,v_0)\widetilde{\alpha}(u_0,v_0)+b_2(u_0,v_0)\widetilde{\beta}(u_0,v_0)+\lambda_v(u_0,v_0)=0,\\
&\widetilde{\alpha}^2(u_0,v_0)+\widetilde{\beta}^2(u_0,v_0)+\widetilde{\gamma}^2(u_0,v_0)=1.
\end{align*}
Since $\mathrm{rank}\mathcal{A}(u_0,v_0)=1,~ a_1^2(u_0,v_0)+b_1^2(u_0,v_0)>\lambda_u^2(u_0,v_0)~ \mathrm{or}~a_2^2(u_0,v_0)+b_2^2(u_0,v_0)>\lambda_v^2(u_0,v_0)$, without loss of generality, we assume $a_1(u_0,v_0)\neq0$ and $a_1^2(u_0,v_0)+b_1^2(u_0,v_0)>\lambda_u^2(u_0,v_0)$. Then the equations
\begin{align*}
&a_1(u_0,v_0)\widetilde{\alpha}(u_0,v_0)+b_1(u_0,v_0)\widetilde{\beta}(u_0,v_0)+\lambda_u(u_0,v_0)=0,\\
&a_2(u_0,v_0)\widetilde{\alpha}(u_0,v_0)+b_2(u_0,v_0)\widetilde{\beta}(u_0,v_0)+\lambda_v(u_0,v_0)=0
\end{align*}
are equivalent to
\begin{align*}
a_1(u_0,v_0)\widetilde{\alpha}(u_0,v_0)+b_1(u_0,v_0)\widetilde{\beta}(u_0,v_0)+\lambda_u(u_0,v_0)=0.
\end{align*}
Since $\widetilde{\alpha}^2(u_0,v_0)+\widetilde{\beta}^2(u_0,v_0)+\widetilde{\gamma}^2(u_0,v_0)=1$, we calculate that
\begin{align*}
\widetilde{\alpha}(u_0,v_0)=&-\frac{b_1\widetilde{\beta}+\lambda_u}{a_1}(u_0,v_0),\\
\widetilde{\gamma}(u_0,v_0)=&\pm\frac{-(a_1^2+b_1^2)\widetilde{\beta}^2
-2\lambda_ub_1\widetilde{\beta}+a_1^2-\lambda_u^2}{a_1}(u_0,v_0).
\end{align*}
Note that $\widetilde{\alpha}^2(u_0,v_0)+\widetilde{\beta}^2(u_0,v_0)\leq1$, we have $$\widetilde{\beta}(u_0,v_0)\in\bigg[\frac{-\lambda_ub_1-|a_1|\sqrt{a_1^2+b_1^2-\lambda_u^2}}{a_1^2+b_1^2}(u_0,v_0),
\frac{-\lambda_ub_1+|a_1|\sqrt{a_1^2+b_1^2-\lambda_u^2}}{a_1^2+b_1^2}(u_0,v_0)\bigg].$$
It follows that
$$
\widetilde{\bm{\nu}}(u_0,v_0)=\widetilde{\alpha}(u_0,v_0) \bm{s}(u_0,v_0)+\widetilde{\beta}(u_0,v_0) \bm{t}(u_0,v_0)+\widetilde{\gamma}(u_0,v_0) \bm{n}(u_0,v_0)
$$
is a circle on the unit sphere. Therefore, $\bm{x}(u_0,v_0)+\lambda(u_0,v_0)\widetilde{\bm{\nu}}(u_0,v_0)$  represents a circle on the sphere $S_{(x(u_0,v_0), \lambda(u_0,v_0))}$ that can be denoted as $C_{(x(u_0,v_0), \lambda(u_0,v_0))}$.
We prove the following theorem which asserts that
$\bm{f}=\mathcal{D}$ if and only if $\bm{x}:U\to \mathbb{R}^3$ is regular,
and $\mathcal{D}$ contains not only $\bm{f}$ but also a circle
or a sphere at a singular point $(u,v)$ of $\bm{x}$
when $\bm{x}$ is singular.
\begin{theorem}\label{theorem3}
Let $(\bm{x},\bm{n},\bm{s}): U\to \mathbb{R}^3\times\Delta$, $\lambda: U\to \mathbb{R}_+$
be a framed surface and a positive function respectively.
Suppose that the sphere family
$\mathcal{S}_{(x, \lambda)}$ is creative.
Then, the following holds.
\[
\mathcal{D}=\bm{f}\cup
\left(\bigcup_{(u,v)\in \Sigma_4}C_{(x(u,v), \lambda(u,v))}\right)\cup
\left(\bigcup_{(u,v)\in \Sigma_5}S_{(x(u,v), \lambda(u,v))}\right).
\]
\end{theorem}
\begin{proof}
Recall that $\mathcal{D}$ is the set
\begin{align*}
\bigg\{(x_1, x_2, x_3)\in \mathbb{R}^3\, |\,
\exists (u,v)\in U \mbox{ s.t. } F(x_1, x_2, x_3, u,v)=\frac{\partial F}{\partial u}(x_1, x_2, x_3, u,v)
=\frac{\partial F}{\partial v}(x_1, x_2, x_3, u,v)=0
\bigg\},
\end{align*}
where $F(x_1, x_2, x_3, u, v)=||(x_1, x_2, x_3)-\bm{x}(u,v)||^2 - \lambda^2(u,v)$.
Let $\left(x_{10}, x_{20}, x_{30}\right)\in\mathcal{D}$.
It follows that there exists a $(u_0,v_0)\in U$
such that the following (a), (b) and (c) are satisfied.
\begin{enumerate}
\item[(a)]
$\left((x_{10}, x_{20}, x_{30})-\bm{x}(u_0,v_0)\right)
\cdot \left((x_{10}, x_{20}, x_{30})-\bm{x}(u_0,v_0)\right)-\left(\lambda(u_0,v_0)\right)^2=0$.
\item[(b)] $(\partial/\partial u)\left(\left((x_{10}, x_{20}, x_{30})-\bm{x}(u_0,v_0)\right)
\cdot \left((x_{10}, x_{20}, x_{30})-\bm{x}(u_0,v_0)\right)-\left(\lambda(u_0,v_0)\right)^2
\right)=0.$
\item[(c)] $(\partial/\partial v)\left(\left((x_{10}, x_{20}, x_{30})-\bm{x}(u_0,v_0)\right)
\cdot \left((x_{10}, x_{20}, x_{30})-\bm{x}(u_0,v_0)\right)-\left(\lambda(u_0,v_0)\right)^2
\right)=0.$
\end{enumerate}
The condition (a) implies that there exists a 
unit vector $\bm{\nu}_1(u_0,v_0)\in S^2$ 
such that
\[
\left(x_{10}, x_{20}, x_{30}\right)=\bm{x}(u_0,v_0)+\lambda(u_0,v_0)\bm{\nu}_1(u_0,v_0).
\]
The conditions (b) and (c) imply 
\begin{align*}
&\bm{x}_u(u_0,v_0)\cdot \left(\left(x_{10}, x_{20}, x_{30}\right)-\bm{x}(u_0,v_0)\right)
+
\lambda_u(u_0,v_0)\lambda(u_0,v_0)=0,\\
&\bm{x}_v(u_0,v_0)\cdot \left(\left(x_{10}, x_{20}, x_{30}\right)-\bm{x}(u_0,v_0)\right)
+
\lambda_v(u_0,v_0)\lambda(u_0,v_0)=0.
\end{align*}
By substituting,
we have the following.
\begin{align*}
&\lambda(u_0,v_0)\big((\bm{x}_u(u_0,v_0)\cdot\bm{\nu}_1(u_0,v_0))+\lambda_u(u_0,v_0)\big)=0,\\
&\lambda(u_0,v_0)\big((\bm{x}_v(u_0,v_0)\cdot\bm{\nu}_1(u_0,v_0))+\lambda_v(u_0,v_0)\big)=0.
\end{align*}
Since $\lambda(u,v)>0$ for any $(u,v)\in U$, it follows
\begin{align*}
\bm{x}_u(u_0,v_0)\cdot\bm{\nu}_1(u_0,v_0)+\lambda_u(u_0,v_0)=0,~~~~~~
\bm{x}_v(u_0,v_0)\cdot\bm{\nu}_1(u_0,v_0)+\lambda_v(u_0,v_0)=0.
\end{align*}
Let $\bm{\nu}_1(u_0,v_0)=\alpha_1(u_0,v_0) \bm{s}(u_0,v_0)+\beta_1(u_0,v_0) \bm{t}(u_0,v_0)+\gamma_1(u_0,v_0) \bm{n}(u_0,v_0)$ where $\alpha_1^2(u_0,v_0)+\beta_1^2(u_0,v_0)+\gamma_1^2(u_0,v_0)=1$, we have
\begin{align*}
&a_1(u_0,v_0)\alpha_1(u_0,v_0)+b_1(u_0,v_0)\beta_1(u_0,v_0)+\lambda_u(u_0,v_0)=0,\\
&a_2(u_0,v_0)\alpha_1(u_0,v_0)+b_2(u_0,v_0)\beta_1(u_0,v_0)+\lambda_v(u_0,v_0)=0.
\end{align*}

On the other hand, since $\mathcal{S}_{(x, \lambda)}$ is creative,
there must exist a smooth unit vector field
$\bm{\nu}: U\to S^2$ along $\bm{x}: U\to \mathbb{R}^3$ such that
\begin{align*}
 \bm{x}_u(u,v)\cdot\bm{\nu}(u,v)+\lambda_u(u,v)=0,~~~~~~
 \bm{x}_v(u,v)\cdot\bm{\nu}(u,v)+\lambda_v(u,v)=0.
\end{align*}
It is equivalent to say that there exists a mapping $(\alpha,\beta,\gamma):U\to\mathbb{R}^3$ such that the following hold for any $(u,v)\in U$.
\begin{align*}
&a_1(u,v)\alpha(u,v)+b_1(u,v)\beta(u,v)+\lambda_u(u,v)=0,\\
&a_2(u,v)\alpha(u,v)+b_2(u,v)\beta(u,v)+\lambda_v(u,v)=0,\\
&\bm{\nu}(u,v)=\alpha(u,v)\bm{s}(u,v)+\beta(u,v)\bm{t}(u,v)+\gamma(u,v)\bm{n}(u,v).
\end{align*}
Suppose that the parameter $(u_0,v_0)\in U$ is a regular point of $\bm{x}$. Then, $J_F(u_0,v_0)\neq0$, we have $\alpha_1(u_0,v_0)=\alpha(u_0,v_0)$, $\beta_1(u_0,v_0)=\beta(u_0,v_0)$, $\gamma_1(u_0,v_0)=\gamma(u_0,v_0)$ and $\bm{\nu}_1(u_0,v_0)=\bm{\nu}(u_0,v_0)$. Therefore, by the assertion (2) of Theorem \ref{theorem1}, it follows
\[
\left(x_{10}, x_{20}, x_{30}\right)\in \bm{f}.
\]
Suppose that the parameter $(u_0,v_0)\in\Sigma_4$ is a singular point of $\bm{x}$. Since $\mathrm{rank}\mathcal{A}(u_0,v_0)=1,~ a_1^2(u_0,v_0)+b_1^2(u_0,v_0)>\lambda_u^2(u_0,v_0)~ \mathrm{or}~a_2^2(u_0,v_0)+b_2^2(u_0,v_0)>\lambda_v^2(u_0,v_0)$, we have $\alpha_1(u_0,v_0)=\widetilde{\alpha}(u_0,v_0)$, $\beta_1(u_0,v_0)=\widetilde{\beta}(u_0,v_0)$, $\gamma_1(u_0,v_0)=\widetilde{\gamma}(u_0,v_0)$ and $\bm{\nu}_1(u_0,v_0)=\widetilde{\bm{\nu}}(u_0,v_0)$.
It follows that $$\mathcal{D}_0=C_{(x(u_0,v_0), \lambda(u_0,v_0))},$$
where $\mathcal{D}_0=\bigg\{(x_1, x_2, x_3)\in \mathbb{R}^3\, |\,
F(x_1, x_2, x_3, u_0,v_0)=\frac{\partial F}{\partial u}(x_1, x_2, x_3, u_0,v_0)
=\frac{\partial F}{\partial v}(x_1, x_2, x_3, u_0,v_0)=0
\bigg\}.$\\
Suppose that the parameter $(u_0,v_0)\in\Sigma_5$ is a singular point of $\bm{x}$. Since $\mathrm{rank}\mathcal{A}(u_0,v_0)=0,~ \lambda_u(u_0,v_0)=0,~\lambda_v(u_0,v_0)=0$, for any unit vector
${\bf v}\in S^2$, the following equations hold.
\begin{align*}
x_u(u_0,v_0)\cdot{\bf v}+\lambda_u(u_0,v_0)=0,~~~~~~
x_v(u_0,v_0)\cdot{\bf v}+\lambda_v(u_0,v_0)=0.
\end{align*}
Hence, we may choose
any unit vector ${\bf v}\in S^2$ as the unit vector $\bm{\nu}_1(u_0,v_0)$. It follows
\[
\mathcal{D}_0= S_{(x(u_0,v_0), \lambda(u_0,v_0))}.
\]

\end{proof}

\section{Applications}\label{section5}
In \cite{NaT}, the second author and N. Nakatsuyama defined evolutes of framed surfaces. As an application of Theorem \ref{theorem1}, we investigate the relation between envelopes of sphere families and  evolutes of framed
surfaces. We denote a framed surface set as
$$\mathcal{F}(U,\mathbb{R}^3\times\Delta):=\{(\bm{x},\bm{n},\bm{s})\in C^{\infty}(U,\mathbb{R}^3\times\Delta)|(\bm{x},\bm{n},\bm{s})~\text{is a framed surface}\}.$$

\begin{definition}\label{caustic}
Let $(\bm{x},\bm{n},\bm{s}):U\rightarrow\mathbb{R}^3\times\Delta$ be a framed surface with the basic invariants $(\mathcal{A}, \mathcal{F}_1, \mathcal{F}_2)$. If there exits a mapping $\mathcal{C}:\mathcal{F}(U,\mathbb{R}^3\times\Delta)\rightarrow\mathcal{F}(U,\mathbb{R}^3\times\Delta)$, $\mathcal{C}(\bm{x},\bm{n},\bm{s})=(\bar{\bm{x}},\bar{\bm{n}},\bar{\bm{s}})$ defined by
\begin{align*}
&\bar{\bm{x}}(u,v)=\bm{x}(u,v)+\delta(u,v)\bm{n}(u,v),\\
&\bar{\bm{n}}(u,v)=\sin\theta(u,v)\bm{s}(u,v)+\cos\theta(u,v)\bm{t}(u,v),\\
&\bar{\bm{s}}(u,v)=\bm{n}(u,v),
\end{align*}
where $\delta,\theta:U\rightarrow\mathbb{R}$ are functions such that
\begin{align*}
\left(\begin{array}{cc}
        a_1(u,v)+\delta(u,v)e_1(u,v) & b_1(u,v)+\delta(u,v)f_1(u,v) \\
        a_2(u,v)+\delta(u,v)e_2(u,v) & b_2(u,v)+\delta(u,v)f_2(u,v)
      \end{array}
\right)\left(\begin{array}{c}
               \sin\theta(u,v) \\
               \cos\theta(u,v)
             \end{array}
\right)=\left(\begin{array}{c}
                0 \\
                0
              \end{array}
\right)
\end{align*}
for any $(u,v)\in U$, then we call $\bar{\bm{x}}:U\rightarrow\mathbb{R}^3$ an {\it  evolute} of the framed surface $(\bm{x},\bm{n},\bm{s})$.
\end{definition}

\begin{theorem}\label{thcaustic}
Let $\mathcal{S}_{(x, \lambda)}$ be a sphere family such that $\Sigma_2$ is dense in $U$ or $\Sigma_3$ is dense in $U$.
Suppose that $\bm{f}$ is a created envelope and $(\bm{f},\bm{\nu},\bm{\omega}):U\rightarrow\mathbb{R}^3\times\Delta$ is a framed surface. Then $\bm{x}$ is an evolute of $(\bm{f},\bm{\nu},\bm{\omega})$.
\end{theorem}
\begin{proof}
By the assumption and the assertion (3-i) of Theorem \ref{theorem1}, it follows $\bm{f}$ is the unique envelope of the sphere family $\mathcal{S}_{(x, \lambda)}$. This implies $\bm{\omega}(u,v)=\bm{n}(u,v)$ and $\bm{n}(u,v)\cdot\bm{\nu}(u,v)=0$ hold for any $(u,v)\in U$. Thus, $(\bm{f},\bm{\nu},\bm{n}):U\rightarrow\mathbb{R}^3\times\Delta$ is framed surface, by Proposition \ref{finvariants}, $\{\bm{\nu}(u,v),\bm{n}(u,v),\bm{\mu}(u,v)\}$ is a moving frame along $\bm{f}(u,v)$, where $\bm{\mu}(u,v)=\bm{\nu}(u,v)\times \bm{n}(u,v)$. We have the following systems of differential equations:
\begin{align*}
\left(
  \begin{array}{c}
    \bm{f}_u \\
    \bm{f}_v \\
  \end{array}
\right)&=\left(
          \begin{array}{cc}
            a_{f1} & b_{f1} \\
            a_{f2} & b_{f2} \\
          \end{array}
        \right)
        \left(
          \begin{array}{c}
            \bm{n} \\
            \bm{\mu} \\
          \end{array}
        \right),
\\
\left(
  \begin{array}{c}
    \bm{\nu}_u \\
    \bm{n}_u \\
    \bm{\mu}_u \\
  \end{array}
\right)&=\left(
          \begin{array}{ccc}
            0 & e_{f1} & f_{f1} \\
            -e_{f1} & 0 & g_{f1} \\
            -f_{f1} & -g_{f1} & 0  \\
          \end{array}
        \right)
        \left(
          \begin{array}{c}
            \bm{\nu} \\
            \bm{n} \\
            \bm{\mu} \\
          \end{array}
        \right),\\
\left(
  \begin{array}{c}
    \bm{\nu}_v \\
    \bm{n}_v \\
    \bm{\mu}_v \\
  \end{array}
\right)&=\left(
          \begin{array}{ccc}
            0 & e_{f2} & f_{f2} \\
            -e_{f2} & 0 & g_{f2} \\
            -f_{f2} & -g_{f2} & 0  \\
          \end{array}
        \right)
        \left(
          \begin{array}{c}
            \bm{\nu} \\
            \bm{n} \\
            \bm{\mu} \\
          \end{array}
        \right),
\end{align*}
where $a_{fi}$, $b_{fi}$, $e_{fi}$, $f_{fi}$, $g_{fi}$ (i=1,2) are the basic invariants of the framed surface $(\bm{f},\bm{\nu},\bm{n})$.  Since $\bm{n}(u,v)\cdot\bm{\nu}(u,v)=0$, $(\bm{x},\bm{n},\bm{\nu}):U\rightarrow\mathbb{R}^3\times\Delta$ is also a framed surface and there exists a constant function $\theta(u,v)=\pi/2$ such that
\begin{align}\label{equ6}
\bm{n}(u,v)=\sin\theta(u,v)\bm{n}(u,v)+\cos\theta(u,v)\bm{\mu}(u,v).
\end{align}
On the other hand, by the assertion (2) of Theorem \ref{theorem1}, we have $\bm{x}(u,v)=\bm{f}(u,v)-\lambda(u,v)\bm{\nu}(u,v)$. It follows that
\begin{align*}
&(\bm{f}_u(u,v)-\lambda_u(u,v)\bm{\nu}(u,v)-\lambda(u,v)\bm{\nu}_u(u,v))\cdot \bm{n}(u,v)=0,\\
&(\bm{f}_v(u,v)-\lambda_v(u,v)\bm{\nu}(u,v)-\lambda(u,v)\bm{\nu}_v(u,v))\cdot \bm{n}(u,v)=0
\end{align*}
from $\bm{x}_u(u,v)\cdot \bm{n}(u,v)=0$ and $\bm{x}_v(u,v)\cdot \bm{n}(u,v)=0.$ By substituting equation (\ref{equ6}), we obtain
\begin{align*}
a_{f1}(u,v)-\lambda(u,v)e_{f1}(u,v)=0,~~a_{f2}(u,v)-\lambda(u,v)e_{f2}(u,v)=0
\end{align*}
hold for any $(u,v)\in U.$ Then, we have
\begin{align*}
&\left(\begin{array}{cc}
        a_{f1}(u,v)-\lambda(u,v)e_{f1}(u,v) & b_{f1}(u,v)-\lambda(u,v)f_{f1}(u,v) \\
        a_{f2}(u,v)-\lambda(u,v)e_{f2}(u,v) & b_{f2}(u,v)-\lambda(u,v)f_{f2}(u,v)
      \end{array}
\right)\left(\begin{array}{c}
               \sin\theta(u,v) \\
               \cos\theta(u,v)
             \end{array}
\right)\\
=&\left(\begin{array}{cc}
        a_{f1}(u,v)-\lambda(u,v)e_{f1}(u,v) & b_{f1}(u,v)-\lambda(u,v)f_{f1}(u,v) \\
        a_{f2}(u,v)-\lambda(u,v)e_{f2}(u,v) & b_{f2}(u,v)-\lambda(u,v)f_{f2}(u,v)
      \end{array}
\right)\left(\begin{array}{c}
               1 \\
               0
             \end{array}
\right)=\left(\begin{array}{c}
                0 \\
                0
              \end{array}
\right).
\end{align*}
Therefore, by Definition \ref{caustic}, $\bm{x}$ is an evolute of $(\bm{f},\bm{\nu},\bm{\omega})$.
\end{proof}

We can also study the relation between envelopes of sphere families and pedal surfaces of framed surfaces by applying Theorem \ref{theorem1}. We define a $l$-pedal surface of the framed surface as follows.
\begin{definition}
Let $(\bm{x},\bm{n},\bm{s}):U\rightarrow\mathbb{R}^3\times\Delta$ be a framed surface and $\bm{l}:U\rightarrow S^2$ a smooth mapping. The {\it $l$-pedal surface} of $\bm{x}$ relative to $\bm{P}$ is defined by
\begin{align*}
l\text{-}Pe_P[x](u,v)=\left[(\bm{x}(u,v)-\bm{P})\cdot \bm{l}(u,v)\right]\bm{l}(u,v),
\end{align*}
where $\bm{P}\in\mathbb{R}^3$ is a fixed point. Moreover, we call $l$-pedal surface a {\it pedal surface} of $\bm{x}$ if $\bm{l}(u,v)=\bm{n}(u,v).$ Namely,
\begin{align*}
Pe_P[x](u,v)=\left[(\bm{x}(u,v)-\bm{P})\cdot \bm{n}(u,v)\right]\bm{n}(u,v).
\end{align*}
\end{definition}

\begin{theorem}\label{thpedal}
Suppose a sphere family $\mathcal{S}_{(x, \lambda)}$ creates two envelopes $\bm{f}_1$ and $\bm{f}_2$. If $\Sigma_1$ is dense in $U$ and $\bm{f}_1$ is a constant, then we have
\begin{align*}
\bm{f}_2(u,v)=Pe_{f_1}[2x-\bm{f}_1](u,v).
\end{align*}
\end{theorem}
\begin{proof}
By the assertion (2) of Theorem  \ref{theorem1}, the envelope $\bm{f}(u,v)$ of $\mathcal{S}_{(x, \lambda)}$ satisfies
\begin{align*}
(\bm{f}(u,v)-\bm{x}(u,v))(\bm{f}(u,v)-\bm{x}(u,v))=\lambda^2(u,v).
\end{align*}
Without loss of generality, we choose $\bm{f}_1$ as the origin. Then, we obtain $\bm{x}(u,v)\cdot \bm{x}(u,v)=\lambda^2(u,v).$ It follows that
\begin{align}\label{e7}
\bm{f}(u,v)\cdot(\bm{f}(u,v)-2\bm{x}(u,v))=0.
\end{align}
Furthermore, since $\bm{f}$ is an envelope of $\mathcal{S}_{(x, \lambda)}$, we have the following equalities.
\begin{align}\label{e8}
\bm{f}_u(u,v)\cdot(\bm{f}(u,v)-\bm{x}(u,v))=0,~~~~~~\bm{f}_v(u,v)\cdot(\bm{f}(u,v)-\bm{x}(u,v))=0.
\end{align}
From equalities (\ref{e8}) and by differentiating (\ref{e7}), we obtain
\begin{align*}
\bm{f}(u,v)\cdot \bm{x}_u(u,v)=0,~~~~~~\bm{f}(u,v)\cdot \bm{x}_v(u,v)=0.
\end{align*}
Since $\Sigma_1$ is dense in $U$, it follows there exists a function $\varphi(u,v)$ such that $\bm{f}(u,v)=\varphi(u,v)\bm{n}(u,v)$. Thus, we have $\varphi(u,v)=0$ or $\varphi(u,v)=2\bm{x}(u,v)\cdot \bm{n}(u,v)$ from $\bm{f}(u,v)\cdot(\bm{f}(u,v)-2\bm{x}(u,v))=0.$ Therefore, we obtain
\begin{align*}
\bm{f}_1(u,v)=\mathbf{0},~~~~~~\bm{f}_2(u,v)=(2\bm{x}(u,v)\cdot \bm{n}(u,v))\bm{n}(u,v)=Pe_{f_1}[2x-\bm{f}_1](u,v).
\end{align*}
\end{proof}

According to Theorems \ref{thcaustic} and \ref{thpedal}, we have the following corollary to describe the relation of evolutes and pedal surfaces of framed surfaces.
\begin{corollary}
Let $(\bm{x},\bm{n},\bm{s}),(\bm{f},\bm{\nu},\bm{\omega}):U\rightarrow\mathbb{R}^3\times\Delta$ be two framed surfaces such that $\bm{x}$ is an evolute of $(\bm{f},\bm{\nu},\bm{\omega})$.  Suppose that the set of regular points of $\bm{x}$ is dense in $U$ and $\bm{x}(u,v)\cdot\bm{n}(u,v)\neq0$ holds for any $(u,v)\in U$. Then we have
\begin{align*}
Pe_{P}[x](u,v)=n\text{-}Pe_{P}[f](u,v),
\end{align*}
where $\bm{P}\notin \bm{x}(u,v)$ is a fixed point.
\end{corollary}
\begin{proof}
By the assumption $\bm{P}\notin \bm{x}(u,v)$, we have that $\|\frac{\bm{x}(u,v)-\bm{P}}{2}\|>0$ always holds. Consider the sphere family $\mathcal{S}_{(\frac{x+P}{2}, \|\frac{x-P}{2}\|)}$, it follows $\bm{P}$ is an envelope $\bm{f}_1$ of $\mathcal{S}_{(\frac{x+P}{2}, \|\frac{x-P}{2}\|)}$ from
\begin{align*}
\left(\bm{P}-\frac{\bm{x}(u,v)+\bm{P}}{2}\right)\cdot\left(\bm{P}-\frac{\bm{x}(u,v)+\bm{P}}{2}\right)=\left\|\frac{\bm{x}(u,v)-\bm{P}}{2}\right\|^2
\end{align*}
holds for any $(u,v)\in U.$ Without loss of generality, we choose $\bm{P}$ as the origin. Since the set of regular points of $\bm{x}$ is dense in $U$ and $\bm{x}(u,v)\cdot\bm{n}(u,v)\neq0$ holds for any $(u,v)\in U$, the set $\Sigma_1$ for the sphere family $\mathcal{S}_{(\frac{x+P}{2}, \|\frac{x-P}{2}\|)}$ is dense in $U$. By Theorem \ref{thpedal}, we have
\begin{align*}
\bm{f}_2(u,v)=(\bm{x}(u,v)\cdot \bm{n}(u,v))\bm{n}(u,v)=Pe_{P}[x](u,v).
\end{align*}
On the other hand, since $\bm{x}$ is an evolute of $(\bm{f},\bm{\nu},\bm{\omega})$ and the set of regular points of $\bm{x}$ is dense in $U$, by Theorem \ref{thcaustic}, it follows that there exists a function $\delta(u,v)$ such that $\bm{x}(u,v)=\bm{f}(u,v)+\delta(u,v)\bm{\nu}(u,v)$ and $\bm{n}(u,v)\cdot\bm{\nu}(u,v)=0.$ Thus, we have
\begin{align*}
\bm{f}_2(u,v)=&Pe_{P}[x](u,v)=(\bm{x}(u,v)\cdot \bm{n}(u,v))\bm{n}(u,v)\\
=&(\bm{f}(u,v)\cdot \bm{n}(u,v))\bm{n}(u,v)=n\text{-}Pe_{P}[f](u,v).
\end{align*}
\end{proof}


\section*{Acknowledgement}
{This work was completed during the visit of the third author to Muroran Institute of Technology. We would like to thank Muroran Institute of Technology for their kind hospitality.
The initial version of Proposition \ref{proposition_equivalence} had a crucial error which was pointed
out by Muneteru Sano at Science Tokyo.     Moreover, he kindly suggested a
slightly weak corrected statement with proof.
Based on his suggestion, the authors improved
the statement of Proposition \ref{proposition_equivalence} and its proof.
We are very grateful to him.

The first author  was supported by JSPS KAKENHI (Grant No. 23K03109). The second author was partially supported by JSPS KAKENHI (Grant No. 24K06728). The third author was supported by Fundamental Research Funds for the Central Universities
(Grant No. 3132026202).

}
%
%

\end{document}